\documentclass[12pt,leqno]{amsart}
\usepackage[dvips]{graphicx}
\usepackage{amsfonts}
\usepackage{amssymb}
\usepackage{amsmath}
\usepackage{amscd}
\usepackage{graphicx}
\def\DATE{\today}
\newtheorem{theorem}{Theorem}
\newtheorem{definition}[theorem]{Definition}
\newtheorem{corollary}[theorem]{Corollary}

\newtheorem{lemma}[theorem]{Lemma}

\newtheorem{proposition}[theorem]{Proposition}

\newcommand\C{\mathbb{C}}
\newcommand\R{\mathbb{R}}

\newcommand\g{\mathfrak{g}}

\newcommand\K{\mathbb{K}}

\newcommand\pf{\noindent{\it Proof. }}
\newcommand\lr{\left\{ \begin{array}{l}}
\def\ds{\displaystyle}
\renewcommand{\mod}{\mathrm{mod}}

\newcommand{\cc}{\mathcal C }

\title[On filiform Lie algebras]{On filiform Lie algebras. Geometric and algebraic studies \protect \footnote{This work has been presented to the  13-th INTERNATIONAL WORKSHOP ON DIFFERENTIAL GEOMETRY AND ITS APPLICATIONS, Romania,   2017}}
\author{Elisabeth Remm}
\address{Universit\'{e} de Haute Alsace, LMIA, 6 rue des Fr\`{e}res Lumi\`{e}re, 68093 Mulhouse}
\email{elisabeth.remm@uha.fr}

\begin{document}

\maketitle

\begin{abstract}
A finite dimensional filiform $\K$-Lie algebra is a nilpotent Lie algebra $\g$ whose nil index is maximal, that is equal to $\dim \g -1.$ We describe necessary and sufficient conditions for a filiform algebra over an algebraically closed field of characteristic $0$ to admit a contact linear form (in odd dimension) or a symplectic structure (in even dimension). If we fix a Vergne's basis, the set of filiform $n$-dimensional Lie algebras is a closed Zariski subset of an affine space generated by the structure constants associated with this fixed basis. Then this subset is an algebraic variety and we describe in small dimensions the algebraic components. 

\end{abstract}

{\bf Conventions.} All Lie algebras considered in this paper will be defined over an algebraically closed fixed field $\K$ of characteristic 0.

\medskip

\noindent{\bf Introduction.}
The problem of classification up to isomorphism is a substantial problem in the study of finite dimensional Lie algebras, even over an algebraically closed field $\K$ of characteristic $0$. This problem has  a solution if we consider in the simple or semisimple algebras, that is, non abelian with no non-trivial ideals or with no non-zero abelian ideals. Indeed, in 1884, Elie Cartan gave the classification of the complex and real simple  finite dimensional Lie algebras. His work is based on the works of Killing. He shows that this classification is reduced to $4$ classes and $5$ exceptional Lie algebras. The Levi decomposition, which was a conjecture of Killing and Cartan and was proved by Eugenio Elia Levi (1905), states that any finite-dimensional real or complex Lie algebra  is the semidirect product of a solvable ideal and a semisimple subalgebra. From this result, the problem of classification is reduced to the classification of solvable Lie algebras and to the problem of representation of semisimple Lie algebras. Thus we are led to classify the solvable Lie algebras. But there are only few results on this topic. We know this classification up to the dimensions $5$ or $6$. Without effective solutions to these problems, since the structure of a solvable Lie algebra is determined by its nilradical, that is the maximal nilpotent ideal, we are interested by the classification of nilpotent Lie algebras. It is one of the firstly aims  of this work. Nowadays, only the classifications of complex or real nilpotent Lie algebras   of dimension lower or equal to $7$ was established. Even this partial result was presented after numerous tries by various authors and in often different approaches. It is surely unrealistic to hope for classifications in bigger dimensions.  Indeed, in  dimension $7$, we find non isomorphic  families of one parameter of nilpotent Lie algebras 
but also a very large number (more than $100$) of not parametrized and not isomorphic Lie algebras. In bigger dimension, this number of Lie algebras has to give the dizziness and will become very difficult to verify (it is enough to see the different tries of classifications in dimension $7$ and the story of the dimension $7$ to realize it already).  Moreover, when we have a given  Lie algebra, it is often difficult to find this algebra in the official list of the classification because it is not the same invariants which are used, and finding the change of basis is a little bit tedious. Thus to pursue the work of classification up an isomorphism in bigger dimension seems  utopian. Some works showed that particular families of nilpotent algebras were parametrized by tensor spaces. This means that it equivalent to classify Lie algebras and arbitrary bilinear maps. Therefore, it does not seem wise to study particular properties of Lie algebras  by starting on existing classifications. For example, the properties that we consider in this paper are the existence of contact forms or of symplectic forms, and also topological and algebraic properties based on the deformation theory and on the rigidity property. Thus the approach is more original. It was introduced in a previous paper concerning the study of the $k$-step nilpotent Lie algebras. We globally study some reduced families which are invariant by isomorphism and which are  closed, that is defined by a finite polynomial system. For this families, we can define an adapted cohomology and then introduce a notion of local rigidity, that is we consider only deformations of elements of a given family which stay in this family. In \cite{GRKegel}, we have for the first time studied some geometrical properties of $k$-step nilpotent Lie algebras by considering  family adapted to the characteristic sequence of a nilpotent Lie algebra (see e.g. \cite{RBreadth} to a presentation of this invariant).

We consider in this work filiform Lie algebras, that is nilpotent Lie algebra whose nil index is maximal, that is $n-1$ if $n$ is the dimension of the Lie algebra. Of course, for this family we know the classification up to the dimension $7$. In \cite{Goze-Khakim2}, classifications are also given for the dimensions $8$ and $9$ for this particular family of nilpotent Lie algebras. But there was no general consensus on these classifications. Thus we begin our work with these dimensions. We are interested in topological properties, as the rigidity, which we study by developing this notion of restricted cohomology. In particular we come back on a result of \cite{AG} concerning the local rigidity. We describe also in dimension $8$ the family (which is neither open nor closed) of symplectic Lie algebras, that is Lie algebras provided with a bilinear symplectic form. In dimension $9$, we show that there are no rigid filiform Lie algebras and we determine all the contact Lie algebras. We show also that any symplectic $8$-filiform Lie algebra is obtained as a quotient of a contact filiform  Lie algebra by its center and we find again the first given description of symplectic Lie algebra. We do the same thing for the dimensions $10$ and $11$, giving the description of the contact $11$-dimensional filiform Lie algebras and consequently the description of the symplectic $10$-dimensional filiform Lie algebras. We remark also, that in these dimensions, all none of the Lie algebras are rigid. For the general dimension, we determine the family of contact $(2p+1)$-dimensional filiform Lie algebras. We expose also a model for this geometrical property, that is a family of contact Lie algebras such that any filiform contact Lie algebra is a deformation of an algebra of this model. Let us recall also that  the reduction of  the polynomial Jacobi system, that is the system of polynomial equations given by the Jacobi conditions, is the fundamental problem.  We don't have  many tools to find the generators of the ideal generated by these equations in the ring $\K[C_{i,j}^k, 1 \leq i < j \leq n, 1 \leq k \leq n]$. In general this ideal $I$ is not equal to $\sqrt{I}$ and the associated affine scheme is not reduced. We described for the family model (parametrized by $(p-1)$ parameters) a process of reduction and we give the associated reduced system. To end this study, we determine from the family of  $(2p+1)$-dimensional contact filiform Lie algebra the family of symplectic $(2p)$-dimensional filiform Lie algebras and we propose a notion of deformation of symplectic Lie algebra based on deformations of the contact Lie algebra which is a the one dimensional central extension .
Let us note also, that in \cite{BGR}, we have studied filiform Lie algebras admitting a $G$-grading, where $G$ is an abelian group.

\section{Generalities on Filiform Lie algebras}
\subsection{The Vergne's basis} 
Let $\g$ be a $n$-dimensional Lie algebra over the field $\K.$ The ascending central series $\{\mathcal{C}_i \g \} $ of $\g$ is defined by 
$$ \mathcal{C}_0\g=\{0\}, \quad \mathcal{C}_i \g=\{ X \in \g \, / \, [X,\g] \subset \mathcal{C}_{i-1}\g \}$$ 
for $i>0$ and  the descending central series $\{\mathcal{C}^i \g \} $ of $\g$ is defined by
$$ \mathcal{C}^0\g=\g , \quad \mathcal{C}^i \g=[\g, \mathcal{C}^{i-1}\g], \ \ i >0.$$ 

\begin{definition}
The $n$-dimensional Lie algebra $\g$ is called filiform if we have $\dim \mathcal{C}_i\g=i$ for $0 \leq i \leq n-2.$ 
\end{definition}
 A filiform Lie algebra is nilpotent and we have 
 $$\cc_i\g = \mathcal{C}^{n-i-1}\g,  \ \  0 \leq i \leq n-1.$$

\begin{proposition}
Let $\g$ be a $(n+1)$-dimensional filiform Lie algebra. There exists a basis $\{ X_0, X_1, \cdots , X_n\}$ called the Vergne basis of $\g$ such that 
$$\left\{
\begin{array}{ l}
      [X_0,X_i]=X_{i+1} , \, 1\leq i \leq n-1,   \lbrack X_1,X_{n-1}\rbrack=0,\\
         \lbrack X_i,X_j]=\ds \sum_{k \geq i+j} C^k_{i,j} X_k , 
\end{array}
  \right.$$  
\end{proposition}
Another characterization of a filiform Lie algebra is given by its characteristic sequence. In fact if $X$ is a vector of the nilpotent Lie algebra $\frak{g}$, the characteristic sequence $c(X)$ of the adjoint operator $ad X$  is the decreasing sequence of the dimensions of the Jordan blocks of the nilpotent operator $ad X$. The characteristic sequence $c(\frak{g})$ of $\frak{g}$ is the following sequence $max\{c(X), X \in \g - \mathcal{C}^1(\frak{g})\}$, the maximum corresponding to the lexicographic order. Any vector $X$ whose characteristic sequence $c(X)$ of $ad X$  is equal to $c(\frak{g})$  is called characteristic vector of  the nilpotent Lie algebra $\g$ (so according to the lexicographic ordering,  we have $c(Y)\leq c(X)$ for any $Y \in \g$ if $X$ is a characteristic vector).
The $n$-dimensional nilpotent Lie algebra $\g$ is filiform if and only if $c(\frak{g})=n-1.$ If $\{ X_0, X_1, \cdots , X_n\}$ is a Vergne basis of a filiform Lie algebra $\g$, the characteristic sequence $c(X_0)$ of the adjoint operator $ad X_0$
is equal to $(n,1)$. An interesting example of $(n+1)$-dimensional filiform Lie algebra is the Lie algebra $L_{n+1}$ often called the model filiform Lie algebra (\cite{GK}) whose Lie bracket $\mu_0$ is 
\begin{equation}\label{Ln}
\left\{
\begin{array}{ l}
      \mu_0(X_0,X_i)=X_{i+1} , \, 1\leq i \leq n-1,   \\
         \mu_0 (X_i,X_j)=0, \, 1 \leq i < j \leq n.
\end{array}
  \right.
  \end{equation}
 
 \subsection{Geometric structure on filiform Lie algebras}
 
 \subsubsection{Contact and symplectic structures}
  
  Let $\g$ be a $(2p)$-dimensional $\K$-Lie algebra.
A symplectic form on $\g$ is a closed 2-form $\theta$, that is satisfying
$$\theta([X,Y],Z)+\theta([Y,Z],X)+\theta([Z,X],Y)=0$$
for any $X,Y,Z \in \g$ and which is also nondegenerate that is
$$\theta^p= \theta \wedge \cdots \wedge \theta\neq 0.$$
A Lie algebra provided with a symplectic form $\theta$ is called a symplectic Lie algebra and denoted by the pair $(\g,\theta)$. There exist filiform Lie algebras without symplectic structures. For example, in \cite{Bou}, on find the list of filiform Lie algebras up the dimension $6$. When the symplectic form is exact, that is if there exists $\omega$ in $\g^*$ the dual space of $\g$ such that $\theta =d\omega$ where $d\omega$ is the bilinear form defined by $d\omega(X,Y)=-\omega([X,Y])$ for any $X,Y \in \g$, the symplectic Lie algebra is called frobeniusian. But, from \cite{GRcontact}, there are no frobeniusian nilpotent Lie algebras. 

Let $\g$ be a $n=2p+1$-dimensional $\K$-Lie algebra.
A contact form on $\g$ is non zero linear form  $\omega$ on $\g$  satisfying
$$\omega \wedge (d\omega)^p \neq 0$$
where $d\omega$ is the bilinear form defined by $d\omega(X,Y)=-\omega([X,Y])$ for any $X,Y \in \g$ and
$(d\omega)^p=d\omega \wedge \cdots \wedge d\omega)$ $p$-times. In case of nilpotent Lie algebras, there is an obstruction to the existence of contact form (\cite{GRcontact}), the center of $\g$ can be of dimension $1$. But the center of any filiform Lie algebra is always of dimension $1$, then this necessary conditions is always satisfied. Let us note that this does'nt imply that any filiform odd dimensional Lie algebras admits a contact form. 

\begin{proposition}\label{symcon}
Let $(\g,\theta)$ be a $2p$-dimensional filiform symplectic Lie algebra. Then the one dimensional central extension $\g_{\theta} =\g \oplus \K Z$ whose bracket is given by
$$
\left\{
\begin{array}{l}
[X,Y]_{\g_\theta}=[X,Y] + \theta(X,Y)Z, \ \forall X,Y \in  \g  \\
\lbrack X,Z \rbrack_{\g_\theta}=0, \  \forall X \in  \g  \\
 \end{array}
\right.
$$
is a $2p+1$-dimensional filiform contact Lie algebra.
\end{proposition}
\pf  Let $\{X_0,X_1,\cdots,X_{2p-1}\}$ be a Vergne basis of $\g$ and let be $\{\omega_0,\cdots,\omega_{2p-1}\}$ is the dual basis. If $\theta$ is a symplectic form on $\g$ then 
$$d\theta(X_0,X_i,X_{2p-1})=0=\theta(X_{i+1},X_{2p-1}), \ \ i=1,\cdots,2p-1.$$
This implies that $\theta=\omega_{2p-1} \wedge (a_0\omega_0 + a_1 \omega_1) + \theta_1$.

Conversely, if $\g$ is a $(2p+1)$-dimensional contact nilpotent Lie algebra, then its center $Z(\g)$ is one dimensional (\cite{GRcontact}) and the factor algebra $\g/ Z(\g)$ is a symplectic $2p$-dimensional nilpotent Lie algebra. If $\g$ is filiform, then $\g/Z(\g)$ is also filiform. 
In \cite{GKM}, one proves that $\g$  admits a contact form iff the linear form $\omega_{2p}$ is a contact form where $\{\omega_0,\cdots,\omega_{2p}\}$ is the dual basis of the Vergne basis of $\g$. We deduce
\begin{proposition} \label{consym}
A $(2p)$-dimensional Lie algebra $(\g,\theta)$ is symplectic iff in the central extension $\g_{\theta}$, the linear form $\omega_{2p}$ is a contact form, where $\{\omega_0,\cdots,\omega_{2p}\}$ is the dual basis of the Vergne basis of $\g_\theta$.
\end{proposition}

 \subsubsection{Complex structures}
 \begin{definition}A  complex structure on an $(2p)$-dimensional $\R$-Lie algebra $\g$ is a linear endomorphism $J$ of $\g$ such that:
\begin{enumerate}
  \item $J^2 = - Id$;
  \item $ [JX,JY] - [X,Y] - J[JX,Y] - J[X,JY] =0, \qquad \forall X, Y \in \g$.
\end{enumerate}
\end{definition}
 Let $\g_\C = \g \otimes _\R \C$ denote the complexification of $\g$, and 
 $$\sigma : \g_\C \rightarrow \g_\C$$
  the corresponding conjugation. The second condition is  equivalent to the splitting
$$\g_\C = \g^{1,0} \oplus \g^{0,1}$$
where $\g^{1,0}$ and $ \g^{0,1}$ are complex Lie subalgebras of $\g_\C$ and $\g^{0,1} = \sigma (\g^{1,0})$. 

\begin{proposition} \cite{GRComplex} There are no filiform Lie algebra admitting a complex structure. 
\end{proposition}

\subsubsection{Affine structures}
\begin{definition}
An affine structure  on a Lie algebra $\g$ is a $\K$-bilinear multiplication, denoted $X \cdot Y$
 which is left-symmetric, that is  
 $$ X \cdot (Y \cdot Z) - (X \cdot Y) \cdot Z = Y \cdot (X \cdot Z) - (Y \cdot X) \cdot Z$$
for all $X,Y,Z \in \g$.and satisfies
$$[X, Y] =X\cdot Y - Y \cdot X$$
where $[X, Y] $ denotes the Lie bracket of $\g$.
 \end{definition}
 The problem, which concerns also the linear representations of Lie algebra (\cite{Burde}), is no completely solved even for filiform Lie algebras. We know that, as soon as the dimension is greater or equal to $10$, there exist filiform Lie algebras without affine structure. However, let us recall this classical result:
 \begin{proposition}
Any symplectic Lie algebra is affine.
\end{proposition}
\pf Let $(\g,\theta)$ be a symplectic Lie algebra. We consider the product $XY$ given by $XY=f(X)Y$
where $f: \g \rightarrow End(\g)$ is defined implicitly by $\theta (f(X)(Y),Z)=-\theta (Y,[X,Z])$
for any $X,Y,Z \in \g$. Since $\theta$ is symplectic, this product $XY$ is well defined and provides $\g$ with an affine structure.

In the case of contact Lie algebras, see for example (\cite{Remmaffine}).

\section{Filiform Lie algebras of dimension $8$}

\subsection{Topological description}
The classification of filiform  Lie algebras over $\K$ of dimension less than or equal to $7$ is well known (\cite{GK}). The aim of this section is to come back to the classification proposed in \cite{AG} and corrected some mistakes of this last paper.  

Let $\g$ be a $8$-dimensional filiform Lie algebra. If we denote by $\mu$ its Lie bracket and $\{X_0,\cdots,X_7\}$ a Vergne basis, then the Jacobi identity implies 
$$\mu (X_0,\mu(X_i,X_j))=\mu(X_i,X_{j+1})+\mu(X_{i+1},X_j).$$

These identities imply that the structure constants $C_{i,j}^k$ for $k<7$ are linear combinations of $C_{ij}^7.$ We deduce 
\begin{proposition}\cite{AG}\label{p1}
Any $8$-dimensional filiform Lie algebra over $\K$ is given in a Vergne basis by 
\begin{equation}
\label{fili8}
\left\{
\begin{array}{ l}
    \mu(X_0, X_i)=X_{i+1},  \ 1 \leq i \leq 6 ,\\
     \mu(X_2, X_5)= a_1 X_7, \\
       \mu(X_1, X_5)= a_1 X_6+a_2 X_7, \\
 \mu(X_3, X_4)= -a_1 X_7 ,\\
 \mu(X_2, X_4)=  a_4 X_7,\\
  \mu(X_1, X_4)= a_1 X_5 +(a_2+a_4) X_6+a_5X_7,\\
\mu(X_2, X_3)=  a_4 X_6+a_6X_7,\\
\mu(X_1, X_3)= a_1 X_4 +(a_2+2a_4) X_5+(a_5+a_6)X_6+a_7X_7,\\
\mu(X_1, X_2)= a_1 X_3 +(a_2+2a_4) X_4+(a_5+a_6)X_5+a_7X_6+a_8X_7,\\
\end{array}
\right.
\end{equation}
with  $a_1(5a_4+2a_2)=0.$
\end{proposition} 

\noindent {\bf Consequence:} Let $\mathcal{L}ie_8$ be the algebraic variety over $\K$ of 
$8^{3}$-uple $( C_{i,j}^k )$ with $0 \leq i <j \leq  7$ and $0 \leq k \leq 7$ satisfying 
$$\left\{
\begin{array}{l }
  C_{i,j}^k=-C_{j,i}^k  \\
     \ds \sum_{l=0}^7 C_{i,j}^l C_{l,k}^s +C_{j,k}^l C_{l,i}^s +C_{k,i}^l C_{l,i}^s=0, \ \forall s=0,\cdots, 7.
\end{array}
\right.$$
A $8$-dimensional Lie algebra with Lie bracket $\mu$ is identified with a point of $\mathcal{L}ie_8$ considering the structural constants of  $\mu$ in a given the basis. An action of the algebraic group $GL(8,\K)$ on $\mathcal{L}ie_8$ corresponds to the changes of basis and the orbit of a Lie algebra is its isomorphism class. Let $\mathcal{N}il_8$ be the algebraic subvariety of $\mathcal{L}ie_8$ whose elements correpond the $8$-dimensional nilpotent Lie algebra and $\mathcal{F}il_8$ be the set of $8$-dimensional filiform Lie algebras. It is a Zariski open set of $\mathcal{N}il_8$ and from Proposition \ref{p1} it is the orbit of the subvariety $Fil_8$ of $\mathcal{N}il_8$  whose elements are the Lie algebras defined in (\ref{fili8}). We deduce that the study of the open set $\mathcal{F}il_8$ can be deduced directly from 
the study of $Fil_8$. 

\medskip

 The set $Fil_8$ is an algebraic variety embedded in $\K^8$ and parametrized by the structural constants $a_1, a_2,a_4,a_5,a_6,a_7, a_8$. It is the union of two irreducible connected algebraic components 
 \begin{enumerate}
  \item $Fil_8(1)$ defined by $a_1=0$ which is a $6$-dimensional plane,
  \item $Fil_8(2)$ defined by $5a_4+2a_2=0$  which is also a  $6$-dimensional plane.
\end{enumerate} 

We deduce that $\mathcal{F}il_8$ is the union of two irreducible algebraic components, $\mathcal{F}il_8(1)$ and $\mathcal{F}il_8(2)$ which are respectively the orbits of $Fil_8(1)$ and $Fil_8(2)$.

\subsection{Deformations}
In the following, we identify the Lie bracket of a $8$-dimensional nilpotent (resp. filiform) Lie algebra with the point $(C_{ij}^k)$ of $\mathcal{N}il_8$ (resp. $\mathcal{F}il_8$) where the $C_{i,j}^k$ are its structural constants related to the given Vergne basis $\{X_0,\cdots,X_7\}$.

\begin{definition}
Let $\mu_0$ be a Lie bracket belonging to $\mathcal{F}il_8$. A deformation $\mu$ in $\mathcal{F}il_8$  of $\mu_0$ is a formal deformation in the Gerstenhaber sense such that $\mu \in \mathcal{F}il_8 \otimes \K[[t]]$.
\end{definition}

From \cite{Goze-Khakim}, any deformation in $\mathcal{F}il_8$ is isomorphic to a linear deformation $\mu_0 +t \psi$ where $\psi$ is a bilinear  skew-symmetric application which is a 2-cocycle for the Chevalley-Eilenberg cohomology of $\mu_0$ and satisfying also the Jacobi identity. Moreover since $\mu_0 +t \psi$ is filiform, $\psi$ is a nilpotent Lie bracket.

The determination of deformations in $\mathcal{F}il_8$ of $\mu_0 \in \mathcal{F}il_8$ reduces  to the study of bilinear skew-symmetric applications $\psi$ such that $\mu_0 +t \psi$  is in $\mathcal{F}il_8$.

\begin{lemma}
Let $\mu_0$ be in $Fil_8$. Then $\mu_0 +t \psi$ is a linear deformation in $Fil_8$ of $\mu_0$ if and only if $\psi$ is given by
\begin{equation}
\label{psi8}
\left\{
\begin{array}{ l}
       \psi(X_2, X_5)= u_1 X_7, \\
       \psi(X_1, X_5)= u_1 X_6+u_2 X_7, \\
\psi(X_3, X_4)= -u_1 X_7 ,\\
\psi(X_2, X_4)= u_4 X_7,\\
 \psi(X_1, X_4)= u_1 X_5 +(u_2+u_4) X_6+u_5X_7,\\
\psi(X_2, X_3)= u_4 X_6+u_6X_7,\\
\psi(X_1, X_3)= u_1 X_4 +(u_2+2u_4) X_5+(u_5+u_6)X_6+u_7X_7,\\
\psi(X_1, X_2)= u_1 X_3 +(u_2+2u_4) X_4+(u_5+u_6)X_5+u_7X_6+u_8X_7,\\
\end{array}
\right.
\end{equation}
with  $u_1(5a_4+2a_2)+a_1(5u_4+2u_2)+tu_1(5u_4+2u_2)=0$ where $(a_1,a_2,a_4,a_5,a_6,a_7,a_8)$ are the parameters of $\mu_0$.
\end{lemma}

\subsection{Study of the component $\mathcal{F}il_8(1)$}

Any Lie algebra in $\mathcal{F}il_8(1)$ is isomorphic to a Lie algebra of $Fil_8(1)$ whose Lie bracket is defined by
\begin{equation}
\label{fili8(1)}
\left\{
\begin{array}{ l}
    \mu_0(X_0, X_i)=X_{i+1},  \ 1 \leq i \leq 6 ,\\
       \mu_0(X_1, X_5)= a_2 X_7, \\
 \mu_0(X_2, X_4)=  a_4 X_7,\\
  \mu_0(X_1, X_4)= (a_2+a_4) X_6+a_5X_7,\\
\mu_0(X_2, X_3)= a_4 X_6+a_6X_7,\\
\mu_0(X_1, X_3)= (a_2+2a_4) X_5+(a_5+a_6)X_6+a_7X_7,\\
\mu_0(X_1, X_2)= (a_2+2a_4) X_4+(a_5+a_6)X_5+a_7X_6+a_8X_7.
\end{array}
\right.
\end{equation}

We have seen that $Fil_8(1)$ is an algebraic subvariety of $\mathcal{N}il_8$.  In  \cite{GRKegel}, we have define for each set 
$k$-$\mathcal{N}il_n$ of $k$-step nilpotent $n$-dimensional Lie algebra a cochain complex whose associated cohomology parametrizes the "internal" deformations, that is deformations of $k$-step nilpotent Lie algebras which are also $k$-step nilpotent. When $k$ is maximal, that is for the filiform case, this cohomology is the Vergne cohomology because any nilpotent deformation of a filiform Lie algebra is always filiform. We consider here the same approach, considering the cohomology adapted to the internal deformations in $Fil_8(1)$ that is deformations of elements of $Fil_8(1)$ which remain in this variety. Since we are only concerned by the second space of this cohomology, we shall describe it. Let $\mu$ be in $Fil_8(1)$. We denote by $Z^2_{CR}(\mu,\mu)$  the space of $2$-cochains, that is bilinear skew-symmetric maps $\psi$ on $\K^8$ with values on $\K^8$  which are defined by
\begin{equation}
\label{psi}
\left\{
\begin{array}{ l}
       \psi(X_1, X_5)= u_2 X_7, \\
 \psi(X_2, X_4)=  u_4 X_7,\\
  \psi(X_1, X_4)= (u_2+u_4) X_6+u_5X_7,\\
\psi(X_2, X_3)= u_4 X_6+u_6X_7,\\
\psi(X_1, X_3)= (u_2+2u_4) X_5+(u_5+u_6)X_6+u_7X_7,\\
\psi(X_1, X_2)= (u_2+2u_4) X_4+(u_5+u_6)X_5+u_7X_6+u_8X_7.
\end{array}
\right.
\end{equation}
If $\partial_\mu$ is the coboundary operator of the Chevalley-Eilenberg complex of $\mu$, then $\partial_\mu(\psi)=0$ and any cochain is closed. Let $B^2_{CR}(\mu,\mu)$ the space of $\partial_\mu (f)$ for $f \in End(\K^8)$ such that $f(X_0)=\sum_{i=0}^7 \alpha_i X_i$ and $f(X_1)=\sum_{i=1}^7 \beta_i X_i$. We have
\begin{equation}
\left\{
\begin{array}{ l}
       \delta f(X_1, X_5)= v_2 X_7, \\
 \delta f(X_2, X_4)=  v_4 X_7,\\
  \delta f(X_1, X_4)= (v_2+v_4) X_6+v_5X_7,\\
\delta f(X_2, X_3)= v_4 X_6+v_6X_7,\\
\delta f(X_1, X_3)= (v_2+2v_4) X_5+(v_5+v_6)X_6+v_7X_7,\\
\delta f(X_1, X_2)= (v_2+2v_4) X_4+(v_5+v_6)X_5+v_7X_6+v_8X_7,
\end{array}
\right.
\end{equation}
with 
\begin{equation} \label{df}
\left\{
\begin{array}{ l}
   v_2=a_2(\beta_1-2\alpha_0),\\
    v_4=a_4(\beta_1-2\alpha_0),\\
     v_5=a_5(\beta_1-3\alpha_0)+\alpha_1(-2a_2^2-5a_2a_4-5a_4^2),\\
      v_6=a_6(\beta_1-3\alpha_0)+\alpha_1(-3a_2a_4-3a_4^2),\\
      v_7=  a_7(\beta_1-4\alpha_0)-2a_4 \beta_3-\alpha_1(a_5+a_6)(5a_2+9a_4), \\
      v_8=  a_8(\beta_1-5\alpha_0)-3a_7\alpha_1(2a_2+3a_4)-2a_6 \beta_3-3a_4 \beta_4 \\
   \qquad \qquad    +3\alpha_3a_4(a_2+2a_4)
      -\alpha_1(a_5+a_6)(3a_5+2a_6).\\   
\end{array}
\right.
\end{equation}
Then $B^2_{CR}(\mu,\mu)$ is a linear subspace of $Z^2_{CR}(\mu,\mu)$ and the quotient space 
$$H^2_{CR}(\mu,\mu)=Z^2_{CR}(\mu,\mu) / B^2_{CR}(\mu,\mu)$$
 parametrizes the deformations in $Fil_8(1)$.  We deduce, using the classical theory of Nijenhuis-Richardson, that a Lie algebra $\g$ with bracket $\mu$ is such that $H^2_{CR}(\mu,\mu)=0$ is rigid in $Fil_8(1)$, that is, its orbit in $Fil_8$ is open. But, since $Fil_8(1)$ is isomorphic to a linear $6$-dimensional space, its associated affine scheme is naturally reduced and the converse is also true: if $\g$ is rigid, then $\dim H^2_{CR}(\mu,\mu)=0$. 
But any cocycle $\psi$ with $u_2 u_4 \neq 0$ can not be cohomologous with a cocycle where $u_2= u_4 = 0.$
We deduce that $\ds \dim H^2_{CR}(\mu,\mu) \geq 1.$  We deduce

\begin{proposition}\label{t2}
No Lie algebra belonging to $Fil_8(1)$ is rigid in $\mathcal{F}il_8$ and also in $\mathcal{N}il_8$ and  in $\mathcal{L}ie_8$ where  $\mathcal{N}il_8$ (respectively $\mathcal{L}ie_8)$ is the algebraic variety of $8$-dimensional nilpotent Lie algebras (respectively the algebraic variety of $8$-dimensional Lie algebras).
\end{proposition}

\medskip

Let us determine the Lie algebras of this component which satisfy $\dim H^2_{CR}(\mu,\mu)=1.$ 
Let $\g$ be such a Lie algebra. From the previous remark, any cocycle $\psi \in Z^2_{CR}(\mu,\mu)$ must be cohomologous to a cocycle with $u_5=u_6=u_7=u_8=0$ and $ a_2$ or $a_4$ not zero. Suppose $a_4 \neq 0.$  The coefficients $\beta_3$ and $\beta_4$  can always be chosen such that $u_7-v_7=0$ and $u_8-v_8=0.$ 
\begin{enumerate}
  \item If $a_2+a_4 \neq 0$ and 
$a_5 \neq 0$ we can choose 
$\alpha_1, \beta_1-3\alpha_0$ and $ \beta_1-2\alpha_0$ such that $u_6-v_6=u_5-v_5=u_2-v_2=0.$ The corresponding Lie algebra satisfies $\dim H^2_{CR}(\mu,\mu)=1.$
  \item If $a_2+a_4 \neq 0,$  $a_5 = 0, a_6 \neq 0$ and $2a_2^2+5a_2a_4+5a_4^2\neq 0$ that is $\ds \frac{a_2}{a_4}\neq \frac{-5\pm i\sqrt{15}}{4}$ (here $\K=\C$) we can choose $\alpha_1, \beta_1-3\alpha_0$ and $ \beta_1-2\alpha_0$ such that $u_6-v_6=u_5-v_5=u_2-v_2=0.$ The corresponding Lie algebra satisfies $\dim H^2_{CR}(\mu,\mu)=1.$
  \item $a_2+a_4 = 0$ and $a_6\neq 0$ we can choose $\alpha_1, \beta_1-3\alpha_0$ and $ \beta_1-2\alpha_0$ such that $u_6-v_6=u_5-v_5=u_2-v_2=0.$ The corresponding Lie algebra satisfies $\dim H^2_{CR}(\mu,\mu)=1.$
\end{enumerate}
In all other cases $\dim H^2_{CR}(\mu,\mu)>1.$ 

To simplify denote by $\mu(a_2,a_4,a_5,a_6,a_7,a_8)$ a Lie bracket of a Lie algebra belonging to $Fil_8(1).$
The previous computations shows that the Lie algebras
$\mu(\alpha, 1,0,1,0,0)$ with $\alpha$ such that $2\alpha^2+5\alpha+5\neq 0$ and $\mu(\alpha, 1,1,0,0,0)$ with $\alpha\neq -1 $ satisfy $\dim H^2_{CR}(\mu,\mu)=1.$ 
From (\cite{AG}), these two Lie algebras are isomorphic. So consider the one-parameter family $\mathcal{T}^1_\alpha $ constituted of 
$\mu(\alpha, 1,1,0,0,0)$ with $\alpha\neq -1 .$  Any deformation in $Fil_8(1)$ of an algebra of this family  belongs to this family. Since $Fil_8(1)$ is a $6$-dimensional plane, a reduced algebraic variety, we deduce that the closure of 
$\mathcal{T}^1_\alpha $ is $Fil_8(1).$

\begin{proposition}
The family $\mathcal{T}^1_\alpha $  
\begin{equation}
\label{fili8(1)}
\left\{
\begin{array}{ l}
    \mu(X_0, X_i)=X_{i+1},  \ 1 \leq i \leq 6 ,\\
       \mu(X_1, X_5)= \alpha X_7, \\
 \mu(X_2, X_4)=   X_7,\\
  \mu(X_1, X_4)= (\alpha+1) X_6+X_7,\\
\mu(X_2, X_3)=  X_6,\\
\mu(X_1, X_3)= (\alpha+2) X_5+X_6,\\
\mu(X_1, X_2)= (\alpha+2) X_4+X_5.
\end{array}
\right.
\end{equation}

is rigid in $Fil_8(1)$ and $\mathcal{F}il_8(1)=\overline{\mathcal{O}(\mathcal{T}^1_\alpha )}.$
\end{proposition}

\noindent{\bf Remarks}
\begin{enumerate}
  \item This generalized notion of rigidity which concerns a one-parameter family of Lie algebras has already been defined in(\cite{Goze-Khakim}.
  \item To compare these results with \cite{AG}, we give the complex classification, up to isomorphism, of elements of $Fil_8$. Recall that such results would be utopic to establish for greater dimensions.
\end{enumerate} 
\begin{proposition}
Let us write $(a_2,a_4,a_5,a_6,a_7,a_8)$ a Lie algebra of $Fil_8(1)$. Then any Lie algebra of $Fil_8(1)$ is isomorphic to one of the following:
$$
\begin{array}{llll }
 (\lambda,1,-1,1,0,0)     &  (\lambda,1,0,0,0,0)   &  (-2,1,1,0,0,0)  &  (1,0,-1,1,\lambda,0)     \\
   (0,0,\lambda, 1,1,0)   & (0,0,\lambda, 1,0,0) &  (\lambda,0,0,0,1,1)& (1,0,0,0,1,0)\\
   (1,0,0,0,0,1) & (1,0,0,0,0,0) & (0,0,1,0,1,0) & (0,0,1,0,0,0) \\
   (0,0,0,0,1,0) & (0,0,0,0,0,1) & (0,0,0,0,0,0).
   \end{array}
   $$
   \end{proposition}

\subsection{Study of the component $\mathcal{F}il_8(2)$}

We consider the Lie algebras of (\ref{fili8}) with  $a_4=-\frac{2}{5}a_2.$ Any Lie algebra of $\mathcal{F}il_8(2)$ is isomorphic to
a Lie algebra of $Fil_8(2)$ with Lie bracket defined by:
\begin{equation}
\label{fili8(2)}
\left\{
\begin{array}{ l}
\ds
    \mu(X_0, X_i)=X_{i+1},  \ 1 \leq i \leq 6 ,\\
     \mu(X_2, X_5)= a_1 X_7, \\
       \mu(X_1, X_5)= a_1 X_6+a_2 X_7, \\
 \mu(X_3, X_4)= -a_1 X_7 ,\\
 \mu(X_2, X_4)= - \ds \frac{2}{5} a_2 X_7,\\
  \mu(X_1, X_4)= a_1 X_5 +\ds \frac{3}{5} a_2 X_6+a_5X_7,\\
\mu(X_2, X_3)= \ds -\frac{2}{5} a_2 X_6+a_6X_7,\\
\mu(X_1, X_3)= a_1X_4 +\ds \frac{1}{5} a_2 X_5+(a_5+a_6)X_6+a_7X_7,\\
\mu(X_1, X_2)= a_1 X_3 +\ds \frac{1}{5} a_2 X_4+(a_5+a_6)X_5+a_7X_6+a_8X_7,\\
\end{array}
\right.
\end{equation}

Then $Fil_8(2)$ is a $6$-dimensional plane parametrized by $(a_1,a_2,a_5,a_6,a_7,a_8).$
We consider similarly to the previous section the linear deformations in $Fil_8(2)$ of the Lie brackets belonging to $Fil_8(2).$  
We denote always by $H^2_{CR}(\mu,\mu)$ the space which parametrizes these deformations. The space of $2$-cocycles 
$Z^2_{CR}(\mu,\mu)$ is constituted of the skew-symmetric bilinear applications $\psi$
 given by 
\begin{equation}
\label{psi2}
\left\{
\begin{array}{ l}
       \psi(X_2, X_5)= u_1 X_7, \\
       \psi(X_1, X_5)= u_1 X_6+u_2 X_7, \\
 \psi(X_3, X_4)= -u_1 X_7 ,\\
 \psi(X_2, X_4)= - \ds \frac{2}{5} u_2 X_7,\\
  \psi(X_1, X_4)= u_1 X_5 +\ds \frac{3}{5} u_2 X_6+u_5X_7,\\
\psi(X_2, X_3)= \ds -\frac{2}{5} u_2 X_6+u_6X_7,\\
\psi(X_1, X_3)= u_1X_4 +\ds \frac{1}{5} u_2 X_5+(u_5+u_6)X_6+u_7X_7,\\
\psi(X_1, X_2)= u_1 X_3 +\ds \frac{1}{5} u_2 X_4+(u_5+u_6)X_5+u_7X_6+u_8X_7.\\

\end{array}
\right.
\end{equation}
Let $f$ be a linear  endomorphism of $\g$ such that $\delta f \in Z^2_{CR}(\mu,\mu).$ If $f(X_0)=\sum_{i=0}^7 \alpha_i X_i$ and $f(X_1)=\sum_{i=1}^7 \beta_i X_i$ then
\begin{equation}
\left\{
\begin{array}{ l}
       \delta f   (X_2, X_5)= v_1 X_7, \\
        \delta f  (X_1, X_5)= v_1 X_6+v_2 X_7, \\
   \delta f  (X_3, X_4)= -v_1 X_7 ,\\
  \delta f  (X_2, X_4)= - \ds \frac{2}{5} v_2 X_7,\\
   \delta f  (X_1, X_4)= v_1 X_5 +\ds \frac{3}{5} v_2 X_6+v_5X_7,\\
  \delta f  (X_2, X_3)= \ds -\frac{2}{5} v_2 X_6+v_6X_7,\\
  \delta f  (X_1, X_3)= v_1X_4 +\ds \frac{1}{5} v_2 X_5+(v_5+v_6)X_6+v_7X_7,\\
 \delta f   (X_1, X_2)= v_1 X_3 +\ds \frac{1}{5} v_2 X_4+(v_5+v_6)X_5+v_7X_6+v_8X_7,
\end{array}
\right.
\end{equation}
with 
\begin{equation} \label{df}
\left\{
\begin{array}{ l}
   v_1=a_1(\beta_1-\alpha_0-\alpha_1a_1),\\
    v_2=a_2(\beta_1-2\alpha_0-3\alpha_1a_1),\\
     v_5=a_5(\beta_1-3\alpha_0-5a_1\alpha_1)-\alpha_1(2a_1a_6+\frac{4}{5}a_2^2)+2\alpha_3a_1^2-2\beta_3a_1,\\
      v_6=a_6(\beta_1-3\alpha_0-2a_1\alpha_1)+\alpha_1(a_1a_5+\frac{18}{25}a_2^2)-2\alpha_3a_1^2+2\beta_3a_1,\\      v_7= a_7(\beta_1-4\alpha_0-5a_1\alpha_1)-\frac{7}{5}\alpha_1a_2(a_5+a_6)-\frac{4}{5}\alpha_3a_1a_2+\frac{4}{5}\beta_3a_2,\\   
      v_8=  a_8(\beta_1-5\alpha_0-5a_1\alpha_1)-\alpha_1(\frac{12}{5}a_2a_7+(a_5+a_6)(3a_5+2a_6))+\\
      \qquad \alpha_3(2a_1(a_5+2a_6)-\frac{6}{25}a_2^2)-\frac{4}{5}\alpha_4a_1a_2+2\alpha_5a_1^2-2\beta_3a_6+\frac{6}{5}\beta_4a_2-2\beta_5a_1,\\  \end{array}
\right.
\end{equation}

\begin{enumerate}
  \item If $a_1=a_2=0$ then $\dim H^2_{CR}(\mu,\mu) \geq 3.$ 
  \item If $a_1=0$ and $a_2 \neq 0,$  then $\dim H^2_{CR}(\mu,\mu)\geq 2.$ 
  \item If $a_2=0,a_1 \neq 0$  then $\dim H^2_{CR}(\mu,\mu)\geq 1 .$ 
   \item Assume now $a_2\neq 0,a_1 \neq 0$,  We will develop this case because it is the part of \cite{AG} which presents a mistake. Let us compute the kernel of the linear system $\{v_i=0\}$. Since $a_1a_2\neq0$ then $v_1=v_2=0$ is equivalent to
   $$\beta_1=-\alpha_1a_1,\  \alpha_0=-2\alpha_1a_1.$$  We can also choose $\beta_5$ to obtain $v_8=0$. Then the system is reduced to
   $$
   \left\{
\begin{array}{ l}
     v_5=-\alpha_1(2a_6a_1+\frac{4}{5}a_2^2)+2\alpha_3a_1^2-2\beta_3a_1,\\
      v_6=(3a_6a_1+a_5a_1+\frac{18}{25}a_2^2)\alpha_1-2\alpha_3a_1^2+2\beta_3a_1,\\     
       v_7= (2a_7a_1-\frac{7}{5}a_2(a_5+a_6))\alpha_1-\frac{4}{5}\alpha_3a_1a_2+\frac{4}{5}\beta_3a_2,\\   
      \end{array}
\right.
$$
The matrix of this system is 
$$
M=
\left(
\begin{array}{ccc}
   -2a_6a_1-\frac{4}{5}a_2^2   &  2a_1^2 & -2a_1\\
    3a_6a_1+a_5a_1+\frac{18}{25}a_2^2  &    -2a_1^2 & 2a_1\\
    2a_7a_1-\frac{7}{5}a_2(a_5+a_6) &-\frac{4}{5}a_1a_2 &+\frac{4}{5}a_2
    \end{array}
\right)
$$
Since this matrix is singular, then $\dim {\text{\rm Ker}}M \geq 1$ and $\dim H^2_{CR}(\mu,\mu) \geq 1.$ Since the affine scheme
 associated with this component is reduced, the Lie algebras with $a_1a_2 \neq 0$ are not rigid. We can now study the conditions to have $\dim H^2_{CR}(\mu,\mu) = 1$, this is equivalent to $\text{\rm rank}(M)=2$ that is
 $$a_1(a_5+a_6)-\ds\frac{2}{25}a_2^2 \neq 0$$
 or
 $$\ds a_1a_2\left(\frac{18}{5}a_5+\frac{26}{5}a_6\right)+\frac{72}{125}a_2^3-4a_1^2a_7  \neq 0.$$ 
In particular, we can take $a_1=a_2=1,a_5=-a_6=t.$ For each value of $t$, the dimension of $H_{2,r}(\g,\g)$ of  the corresponding Lie algebra is equal to $1$. We deduce
 \end{enumerate} 
\begin{proposition}
None of the  Lie algebras of $Fil_8(2)$ is rigid in $\mathcal{F}il_8$. This component is the closure of the one dimensional rigid family $\mathcal{T}^2_t(8)$ of Lie algebras isomorphic to
\begin{equation}
\label{fili8rigid}
\left\{
\begin{array}{ l}
   \mu(X_0, X_i)=X_{i+1},  \ 1 \leq i \leq 6 ,\\
     \mu(X_2, X_5)=  X_7, \\
       \mu(X_1, X_5)=  X_6+ X_7, \\
 \mu(X_3, X_4)= -X_7 ,\\
 \mu(X_2, X_4)= - \ds \frac{2}{5} X_7,\\
  \mu(X_1, X_4)= X_5 +\ds \frac{3}{5} X_6+tX_7,\\
\mu(X_2, X_3)= \ds -\frac{2}{5}  X_6-tX_7,\\
\mu(X_1, X_3)= X_4 +\ds \frac{1}{5}  X_5,\\
\mu(X_1, X_2)=  X_3 +\ds \frac{1}{5}  X_4.\\
\end{array}
\right.
\end{equation}
\end{proposition}

\medskip

\noindent {\bf Remark.} The problem of classification of $8$-dimensional filiform Lie algebras has been already solved. We can find this classification in \cite{Goze-Khakim}. A lot of the results previously obtained of course are direct consequence of this classification. But to obtain a general result of the  classification problem is certainly  utopian.  This implies to develop another way. We began a new approach in \cite{GRKegel} by considering subfamilies of $k$-step nilpotent Lie algebras and defining an adapted cohomology of deformations. The previous calculus are made in this way. Nevertheless, since this classification is known and since we have given this classification for the algebras of the first component, it would be  surprising not to give it for the second component.

\begin{proposition}
Let $\g$ be a filiform Lie algebra belonging to
 $\mathcal{F}il_8(2).$ Then any Lie algebra of $Fil_8(2)$ is isomorphic to one of the following corresponding to
 $$(a_1,a_2,a_5,a_6,a_7,a_8) \in \{(1,0,0,0,0,0),(1,0,0,0,1,0) ,(1,0,1,0,\lambda,0),(1,1,\lambda,-2,0,0)\}.$$  
\end{proposition}

 \subsection{Symplectic structures}

We  determine all filiform $8$-dimensional symplectic  Lie algebras.  A direct approach consists to write the conditions related to the existence of a symplectic form for Lie algebras belonging to each component. 
Let $\g \in Fil_8$ and $\theta$ be a closed $2$-form on $\g$.  Let $\{X_0,\cdots,X_7\}$ the Vergne basis (\ref{fili8}). Computing $d\theta(X_0,X_i,X_7)=0$, we obtain that $\theta(X_i,X_7)=0$ for $i=2,\cdots,6.$ As consequence, $\theta(X_i,X_j)=0$ as soon as $i+j \geq 9$ and $d\theta (X_i,X_j,X_k)=0$ is always satisfies for $i+j+k \geq 9$. We call the weight of the equation $d\theta (X_i,X_j,X_k)=0$ the integer $p=i+j+k$ and we solve this equation for $p=8,7,\cdots,3.$ 
\begin{enumerate}
  \item $p=8$. 
  $$
  \left\{
  \begin{array}{ll}
 0= &(3a_1+2a_3)\theta(X_3,X_5)-a_1\theta(X_1,X_7)+a_1\theta(X_2,X_6), \\
  0= &  (2a_1+a_3)\theta(X_3,X_5)-a_3\theta(X_1,X_7). \\  
\end{array}
\right.
$$
  \item $p=7.$
 $$
 \left\{
  \begin{array}{ll}
  0= & \theta(X_2,X_6)+\theta(X_1,X_7),\\
   0= & \theta(X_3,X_5)-a_1\theta(X_0,X_7)+\theta(X_2,X_6),\\
 0=& a_3\theta(X_0,X_7)-\theta(X_3,X_5),\\
 0=& (3a_1+2a_3)\theta(X_3,X_4)-(a_1+a_3)\theta(X_1,X_6)-a_4\theta(X_1,X_7)\\
 & +(2a_1+a_3)\theta(X_2,X_5)+(a_2+a_4)\theta(X_2,X_6). \\
\end{array}
\right.
$$ 
\end{enumerate} 
Then, we have to consider the matrix
$$
M=\left(
\begin{array}{cccc}
 0 &-a_1& a_1   & 3a_1+2a_3   \\
   0   & -a_3&0& 2a_1+a_3  \\
        0    &1&1& 0  \\
              -a_1    &0&1&1   \\
         a_3   &0&0&-1   \\           
                              \end{array}
                              \right)
                              =\left(
\begin{array}{cccc}
 0 &-a_1& a_1   & a_1   \\
   0   & a_1&0& a_1  \\
        0    &1&1& 0  \\
              -a_1    &0&1&1   \\
         -a_1   &0&0&-1   \\           
                              \end{array}
                              \right)
$$
If $\theta$ is symplectic, one of the scalar $\theta(X_0,X_7)$ or  $\theta(X_1,X_7)$ is non zero which is equivalent to say that the rank of $M$ is less than $4$. But $rank M= 4$ if and only if 
$a_1 \neq 0$ so any $8$-dimensional filiform symplectic Lie algebra $\g$ belongs to $Fil_8(1)$. Moreover, the symplectic form satisfies $\theta(X_i,X_j)=0$ for $i+j \geq 8.$ If we compute the relations of weight $6$ we obtain
$$
 \left\{
  \begin{array}{ll}
  0= & \theta(X_2,X_5)+\theta(X_1,X_6)-a_2\theta(X_0,X_7),\\
   0= & \theta(X_3,X_4)+\theta(X_2,X_5)-a_4\theta(X_0,X_7),\\
 0=& (a_2+2a_4)\theta(X_3,X_4)-(a_2+2a_4)\theta(X_2,X_5)+a_4\theta(X_1,X_6). \\
\end{array}
\right.
$$ 
 But $\theta$ is non degenerate if and only if $\theta(0,7)\theta(1,6)\theta(2,5)\theta(3,4) \neq 0.$ This implies
 \begin{enumerate}
  \item $2a_2+5a_4 \neq 0$, that is $\g \in Fil_8(1)$ and $\g \notin Fil_8(2)$ and $a_4(a_2+2a_4)(2a_2-a_4)\neq 0$
  \item  or $a_2=a_4=0.$
\end{enumerate} 
 \begin{proposition}\label{symp8}
 An $8$-dimensional filiform Lie algebra  is symplectic if and only if it is isomorphic to a Lie algebra $\g \in Fil_8(1)-Fil_8(2)$ and  $a_4(a_2+a_4)(2a_2-a_4)(a_2+2a_4)\neq 0$ or $\g \in Fil_8(1)\cap Fil_8(2)$ and  $a_2=a_4= 0.$
 \end{proposition}
 We can note  that any Lie algebra of the rigid family $\mathcal{T}^1_\alpha$ is symplectic, except for three values of $\alpha$ which are $-2,\frac{1}{2},-\frac{5}{2}.$
 
 \medskip
 
 \subsection{Determination of symplectic $8$-dimensional filiform Lie algebras using contact $9$-dimensional filiform Lie algebras.} 
 
 If $(\g,\theta)$ is a $2p$-dimensional symplectic Lie algebra, then the Lie algebra $\g_{\theta}$ the one dimensional  central extension
 $$\g_{\theta}=\g \oplus_{\theta} \K Z.$$
 Recall that the Lie  bracket $mu_1$of $\g_{\theta}$ is given by
 $$\left\{
 \begin{array}{l}
  \mu_1(X,Y)=\theta(X,Y)Z+\mu(X,Y)    \\
    \mu_1(X,Z)=0  
\end{array}
\right.
$$
for any $X,Y \in \g$ and  $\g_{\theta}$ is a contact $(2p+1)$-dimensional Lie algebra and $Z$ generates the center. From Proposition \ref{consym}, $\g$ is a factor algebra $\g_1/Z(\g_1)$  of a filiform contact algebra and the linear form $\omega_{2p}$ is a contact form. 
Then we can determine all the symplectic filiform algebra in dimension $8$ starting from the contact filiform $9$-dimensional Lie algebras. This study is the aim of the last section, but we can already use these results, all the proofs are given in the following section.

\begin{proposition}
Any $9$-dimensional filiform Lie algebra provided with a contact form is isomorphic to a Lie algebra of the following family:
$$
\left\{
\begin{array}{ l}
    \mu(X_0, X_i)=X_{i+1},  \ 1 \leq i \leq 7 ,\\
       \mu(X_1, X_6)= a_2 X_8, \\
 \mu(X_2, X_5)= a_4 X_8,\\
  \mu(X_1, X_5)= (a_2+a_4) X_7+a_5X_8,\\
  
  \mu(X_3, X_4)= a_6X_8\\
  \mu(X_2, X_4)=  (a_4+a_6) X_7+a_7X_8,\\
\mu(X_1, X_4)= (a_2+2a_4+a_6) X_6+(a_5+a_7)X_7+a_8X_8,\\

\mu(X_2, X_3)= (a_4+a_6) X_6+a_7X_7+a_9X_8,\\
\mu(X_1, X_3)= (a_2+3a_4+2a_6) X_5+(a_5+2a_7)X_6+(a_8+a_9)X_7+a_{10}X_8,\\
\mu(X_1, X_2)= (a_2+3a_4+2a_6) X_4+(a_5+2a_7)X_5+(a_8+a_9)X_6+a_{10}X_7+a_{11}X_8,\\
\end{array}
\right.
$$
with $-3a_4^2+2a_6^2+2a_2a_6+a_4a_6=0$ and $a_2a_4a_6 \neq 0.$
\end{proposition}
Since the center is generated by $X_8$, we deduce 
\begin{proposition}
Any symplectic $8$-dimensional filiform Lie algebra is isomorphic to a Lie algebra of the following family
 $$
\left\{
\begin{array}{ l}
    \mu(X_0, X_i)=X_{i+1},  \ 1 \leq i \leq 6 ,\\
  \mu(X_1, X_5)= (a_2+a_4) X_7,\\
  \mu(X_2, X_4)=  (a_4+a_6) X_7,\\
\mu(X_1, X_4)= (a_2+2a_4+a_6) X_6+(a_5+a_7)X_7,\\
\mu(X_2, X_3)= (a_4+a_6) X_6+a_7X_7,\\
\mu(X_1, X_3)= (a_2+3a_4+2a_6) X_5+(a_5+2a_7)X_6+(a_8+a_9)X_7,\\
\mu(X_1, X_2)= (a_2+3a_4+2a_6) X_4+(a_5+2a_7)X_5+(a_8+a_9)X_6+a_{10}X_7\\
\end{array}
\right.
$$
with  $-3a_4^2+2a_6^2+2a_2a_6+a_4a_6=0$ and $a_2a_4a_6 \neq 0.$ This is equivalent to say
\end{proposition}
We come back to the notations (\ref{fili8}) and we put 
$$b_2=a_2+a_4,\ b_4=a_4+a_6, b_5=a_5+a_7, \  b_6=a_7,b_7=a_8+a_9,b_8=a_{10}.$$
Then the conditions $-3a_4^2+2a_6^2+2a_2a_6+a_4a_6=0$ and $a_2a_4a_6 \neq 0$ implies that
$$b_2=b_4=0, \ \text{\rm or} \ b_4(b_2+b_4)(2b_2-b_4)(b_2+2b_4) \neq 0.$$
 
Any symplectic $8$-dimensional filiform Lie algebra is isomorphic to a Lie algebra of the following family
$$
\left\{
\begin{array}{ l}
    \mu(X_0, X_i)=X_{i+1},  \ 1 \leq i \leq 6 ,\\
  \mu(X_1, X_5)= b_2X_7,\\
  \mu(X_2, X_4)=  b_4 X_7,\\
\mu(X_1, X_4)= (b_2+b_4) X_6+b_5X_7,\\
\mu(X_2, X_3)= b_4X_6+b_6X_7,\\
\mu(X_1, X_3)= (b_2+2b_4) X_5+(b_5+b_6)X_6+b_7X_7,\\
\mu(X_1, X_2)= (b_2+2b_4) X_4+(b_5+b_6)X_5+b_7X_6+b_8X_7\\
\end{array}
\right.
$$
with  $b_2=b_4=0$ or $b_4 (b_2+b_4)(b_2-b_4)(b_2+2b_4)\neq 0.$
 We find again all the conditions of Proposition \ref{symp8}.
 
 \medskip
 
 \noindent This last way to determine the symplectic structures permits also the introduce a notion of symplectic deformation. Recall that a deformation of a symplectic Lie algebra can be non symplectic. The simplest example is given by the even dimensional abelian Lie algebra. This algebra is symplectic and any Lie algebra is isomorphic to a deformation of this abelian algebra and it is clear that exists non symplectic Lie algebras as soon as the dimension is strictly greater than $2$. Likewise if a symplectic Lie algebra $\g_1$ s a deformation of a Lie algebra $\g_0$, this last is not necessarily symplectic. Then the classical notion of deformation is not  well adapted to the notion of symplectic structures. But it is not the case for the contact structures. Any deformation of a contact Lie algebra is always a contact Lie algebra (see also \cite{GRcontact}). This remark leads to introduce a restricted notion of deformation that we shall call {\bf symplectic deformation}:
 
 \begin{definition}
 Let $\g_0$ and $\g_1$ be two symplectic $8$-dimensional filiform Lie algebras and let $\g_0^9$ and $\g_1^9$ be $9$-dimensional contact filiform Lie algebras such that $\g_i=\pi (\g^9_i)$, $i=1,2$ where $\pi$ is the canonical projection $\pi:\g \rightarrow \g/Z(\g)$. We shall call that $\g_1$ is a symplectic deformation of $\g_0$ if $\g_1^9$ is a (classical) deformation of $\g_0^9$.
 \end{definition}

 \subsection{Affine structures}

From \cite{Burde}, we know that any $8$-filiform Lie algebra admits an affine structure. To prove this, we construct  affine structure of adjoint type, that is, if $L_i$ denote the linear map $L_i(X)=X_iX$ then $L_0=ad X_0$. Since $L_I$ for $i \geq 2$ is given by $L_i=[L_0,L_{i-1}]$, such affine structure is completely determinate by $L_1$. For exemple if we consider the rigid family $\mathcal{T}^2_t(8)$ in $Fil_8(2)$, we consider for $L_1$ the linear map whose matrix in the Vergne's basis is
$$
\left(
\begin{array}{cccccccc}
      0 &  \alpha_1    & 0      &  \alpha_2       & 0 & 0&0 &0\\
     0&   0  & 0       & 0 & 0 & 0 &0 & 0\\ 
     0&  0    &  0   &0  &0 &0 &0& 0\\
      0&   -\frac{1}{5}  &  0 & 0 & 0  & 0 & 0 & 0\\
      0&   0  & 0       & 0 & 0 & 0 &0 & 0\\ 
      0&    \frac{2(70\alpha_6-25t-42)}{375}  &  0   &0  &0 &0 &0& 0\\
      0&    2\alpha_3  &    \frac{70\alpha_6-25t-42}{375}    & -\frac{2}{25}(5\alpha_6-3)  & \frac{1}{5} &0& 0& 0\\
       0 &   \alpha_4  &   \alpha_3   & \alpha_5  &\frac{t}{2}-\frac{3}{25}(5\alpha_6-3)  &\alpha_6& -\frac{1}{2} &0 \\
\end{array}
\right)
$$
or
$$
\left(
\begin{array}{cccccccc}
      0 &  \alpha_1    & 0      &  \alpha_2       & 0 & 0&0 &0\\
     0&   0  & 0       & 0 & 0 & 0 &0 & 0\\ 
     0&  0    &  0   &0  &0 &0 &0& 0\\
      0&   \frac{3}{5}  &  0 & 0 & 0  & 0 & 0 & 0\\
      0&   0  & \frac{2}{5}       & 0 & 0 & 0 &0 & 0\\ 
      0&    \frac{2(-210\alpha_6+125t-42)}{375}  &  0   &-\frac{2}{5}   &0 &0 &0& 0\\
      0&     \frac{2(-210\alpha_6+125t-42)}{375}  &    \frac{-210\alpha_6+125t-42}{375}  & -\frac{14}{25}(5\alpha_6+1)  & -\frac{3}{5} &0& 0& 0\\
       0 &   \alpha_4  &   \alpha_3   & \alpha_5  &\frac{t}{2}-\frac{21}{25}(5\alpha_6+1)  &\alpha_6& -\frac{1}{2} &0 \\
\end{array}
\right)
$$
Note that all these affine structure are complete, that is the linear map $R_Y: X \rightarrow X \cdot Y$ is  nilpotent for any $Y$, or equivalently the trace of $R_Y$ is zero. Note also that the linear representation of $\g$:
$$\rho: \g \rightarrow End(\g)$$
given by $\rho(X)=L_X$ is not faithful because in all the previous cases $L_{X_7}=0.$ 
\medskip

\noindent{\bf Remark.} Let us consider the polarization of the product $\nabla (X,Y)= X \cdot Y$. We put
$$
\mu(X,Y)=\nabla(X,Y)-\nabla(Y,X), \ \ s(X,Y)=\nabla(X,Y)+\nabla(Y,X).$$
Since $\K$ is of characteristic not $2$, then
$$\nabla(X,Y)= \ds \frac{s(X,Y)+\mu(X,Y)}{2}$$
and  $\mu$ is a Lie bracket. The applications $\mu$ and $s$ are related by the affine condition
$$
  \begin{array}{l}
  A(X,Y,Z)=   \mu(\mu(X,Y),Z)  +\mu(s(Y,Z),X) -\mu(s(Z,X),Y)+2s(\mu(X,Y),Z)\\
     -s(\mu(Y,Z),X) +s(\mu(X,Z),Y)-s(s(Y,Z),X)+s(s(X,Z),Y) =0.
\end{array}
$$
We have also
$$
A(X,Y,Z)+A(Y,Z,X)+A(Z,X,Y)=2(\mu(s(X,Y),Z)+\mu(s(Y,Z),X)+\mu(s(Z,X),Y)).
$$

 \section{Filiform Lie algebras of dimension $9$}

\subsection{The variety $\mathcal{F}il_9$}
Using a similar approach as in dimension $8$, we obtain
\begin{proposition}
Any $9$-dimensional filiform Lie algebra over $\K$ is given in a Vergne basis by 
\begin{equation}
\label{fili9}
\left\{
\begin{array}{ l}
    \mu(X_0, X_i)=X_{i+1},  \ 1 \leq i \leq 7 ,\\
       \mu(X_1, X_6)= a_2 X_8, \\
 \mu(X_2, X_5)= a_4 X_8,\\
  \mu(X_1, X_5)= (a_2+a_4) X_7+a_5X_8,\\
  
  \mu(X_3, X_4)= a_6X_8\\
  \mu(X_2, X_4)=  (a_4+a_6) X_7+a_7X_8,\\
\mu(X_1, X_4)= (a_2+2a_4+a_6) X_6+(a_5+a_7)X_7+a_8X_8,\\

\mu(X_2, X_3)= (a_4+a_6) X_6+a_7X_7+a_9X_8,\\
\mu(X_1, X_3)= (a_2+3a_4+2a_6) X_5+(a_5+2a_7)X_6+(a_8+a_9)X_7+a_{10}X_8,\\
\mu(X_1, X_2)= (a_2+3a_4+2a_6) X_4+(a_5+2a_7)X_5+(a_8+a_9)X_6+a_{10}X_7+a_{11}X_8,\\
\end{array}
\right.
\end{equation}
with $-3a_4^2+2a_6^2+2a_2a_6+a_4a_6=0.$ 
\end{proposition} 

We denote by $Fil_9$ the set of Lie algebras described above. It is clear that $ \mathcal{F}il_9,$ the open set of $9$-dimensional filiform Lie algebras, is the orbit of $Fil_9$ in $\mathcal{N}ilp_9$ associated with the action of the linear group $GL(9, \K).$ This reduces the study of $\mathcal{F}il_9$ to $Fil_9$.
Let $V^{9}$ be the $9$-dimensional $\K$-vector space characterized by the structure constants $\{a_i\}_{2 \leq i \leq 11}, i \neq 3.$

\begin{proposition}
$Fil_9$ is a $8$-dimensional  irreducible algebraic subvariety of $V^{9}$.
\end{proposition} 
 This algebraic variety $Fil_9$ has only one singular point corresponding to $a_i=0$ for all $i$. In all the other points, the tangent space $T_{\mu}(Fil_9)$ to $Fil_9$  is identified to the vector space of $2$-cocycles:
\begin{equation}
\label{ts9}
\left\{
\begin{array}{ l}
    \varphi(X_0, X_i)=0,  \ 1 \leq i \leq 7 ,\\
      \varphi(X_1, X_6)= u_2 X_8, \\
\varphi(X_2, X_5)= u_4 X_8,\\
 \varphi(X_1, X_5)= (u_2+u_4) X_7+u_5X_8,\\
    \varphi(X_3, X_4)= u_6X_8,\\
 \varphi(X_2, X_4)=  (u_4+u_6) X_7+u_7X_8,\\
\varphi(X_1, X_4)= (u_2+2u_4+u_6) X_6+(u_5+u_7)X_7+u_8X_8,\\
\varphi(X_2, X_3)= (u_4+u_6) X_6+u_7X_7+u_9X_8,\\
\varphi(X_1, X_3)= (u_2+3u_4+2u_6) X_5+(u_5+2u_7)X_6+(u_8+u_9)X_7+u_{10}X_8,\\
\varphi(X_1, X_2)= (u_2+3u_4+2u_6) X_4+(u_5+2u_7)X_5+(u_8+u_9)X_6+u_{10}X_7+u_{11}X_8,\\
\end{array}
\right.
\end{equation}
with $u_2a_6+u_4(\frac{1}{2}a_6-3a_4)+u_6(2a_6+a_2+\frac{1}{2}a_4)=0.$
Let $f \in gl(9,\K)$ be an endomorphism. We  put $f(X_0)=\sum \alpha_iX_i$ and $f(X_1)=\sum\beta_iX_i$. Assume that $\delta f(X_0,X_i)=0$. Then $f(X_{i+1})=\mu(f(X_0),X_i)+\mu(X_0,f(X_i))$ for $i=1,\cdots,7.$ The other components of $\delta f$ are
$$
\left\{
\begin{array}{ll}
v_2=&a_2(\beta_1-2\alpha_0), \ v_4=a_4(\beta_1-2\alpha_0), \ v_6=a_6(\beta_1-2\alpha_0), \\
v_5=&a_5(\beta_1-3\alpha_0)-\alpha_1(2a_2^2+9a_4^2+6a_2a_4+5a_4a_6), \\
  v_7=&a_7(\beta_1-3\alpha_0)-\alpha_1(7a_6^2+3a_2a_4+7a_2a_6+11a_4a_6),\\
v_8=&a_8(\beta_1-4\alpha_0)-\alpha_1((5a_2+11a_4+5a_6)a_5+(6a_2+19a_4+10a_6)a_7)-2\beta_3a_4,\\
v_9=&a_9(\beta_1-4\alpha_0)-\alpha_1((3a_4+4a_6)a_5+(4a_2+9a_4+8a_6)a_7)-2\beta_3a_6,\\
v_{10}=&a_{10}(\beta_1-5\alpha_0)-\alpha_1(P_1^1(a_8,a_9)P_2^1(a_2,a_4,a_6)+P_3^2(a_5,a_7))+\alpha_3P_4^2(a_2,a_4,a_6)\\
&-3\beta_4(a_4+a_6)-2\beta_3a_7,\\
v_{11}=&a_{11}(\beta_1-6\alpha_0)-\alpha_1P_5^2(a_2,a_4,a_5,a_6,a_7,a_8,a_9,a_{10})+\alpha_3P_6^2(a_2,a_4,a_5,a_6,a_7)\\
&\alpha_4P_7^2(a_2,a_4,a_6)-2\beta_3a_9-3\beta_4a_7-2\beta_5(2a_4+a_6),\\
\end{array}
\right.
$$
where $P_i^k$ are an homogeneous polynomials of degree $k$.  If we denote by $H^*_{CR}$ the restricted Chevalley cohomology of Lie algebras belonging to $Fil_9$, we have $\dim H^2_{CR}(\mu,\mu)  \geq 1$ for any $\mu \in Fil_9$. We deduce
\begin{proposition} None  of the  $9$-dimensional filiform $\K$-Lie algebras is rigid in $\mathcal{F}il_9$ and also in $\mathcal{N}ilp_9$, and also in $\mathcal{L}ie_9.$
\end{proposition}
\pf In fact for any $\mu \in \mathcal{F}il_9$, $\dim H^2Ð{CR}(\mu, \mu)\neq 0$. Since the affine schema associated with $Fil_9$ is reduced, any rigid Lie algebra in this variety has a trivial cohomology. Then no $\mu \in Fil_9$ is rigid. 

Now we determinate the Lie algebras $\mu$ such that $\dim H^2_{CR}(\mu, \mu)$ is the smallest one. Assume that $a_6(a_4+a_6)(2a_4+a_6)\neq 0.$ Then we can find $\delta f$ such that 
$$u_2+v_2=u_5+v_5=u_9+v_9=u_{10}+v_{10}=u_{11}+v_{11}=0.$$
In this case, the parameters $\alpha_0,\alpha_1,\beta_3,\beta_4,\beta_5$ are fixed and the others relations $u_i+v_i$ cannot be reduced to $0$. Then there exists a represent ant $\varphi$ in the cohomological classe with $u_2=u_5=u_9=u_{10}=u_{11}=0.$ Since $u_2a_6+u_4(\frac{1}{2}a_6-3a_4)+u_6(2a_6+a_2+\frac{1}{2}a_4)=0$, for a such Lie algebra we have $\dim 
H^2_{CR}(\mu,\mu) =2$ and it is the lower bound. We deduce
\begin{theorem}
The variety $\mathcal{F}il_9$ is the closure of the orbit of a rigid $2$-parameters family.
\end{theorem}
Let us  consider the family $\mathcal{T}^9_t$ of Lie algebras $\mu$ defined by $(a_2=\frac{3t^2-t-2}{2},a_4=t,a_5=1,a_6=1,a_8=u,a_9=0,a_{10}=0,a_{11}=0)$ with $t \neq 0,-1,-\frac{1}{2}.$ This family answers to this theorem.

\subsection{Contact $9$-dimensional filiform Lie algebras}

Let $\g$ be a  $9$-dimensional filiform Lie algebra. Let $\{\omega_0,\cdots,\omega_8\}$ the dual basis of $\{X_0,\cdots,X_8\}$. From Proposition \ref{consym}, $\g$ is a contact Lie algebra if and only if $\omega_8$ is a contact form. This is equivalent to 
$$a_2a_4a_6 \neq 0$$
where $a_i$ are the constant structures of $\g$ described in (\ref{fili9}).
We deduce:
\begin{proposition}\label{contact9}
A $9$-dimensional filiform Lie algebra admit a contact form if and only it is isomorphic to a Lie algebra of $Fil_9$ whose structure constant  given in  (\ref{fili9}) satisfy $a_2a_4a_6 \neq 0.$
\end{proposition}

\medskip

\subsection{Come back on $8$-dimensional symplectic filiform Lie algebras.}
The previous theorem is the result announced in Proposition 18. As we have say, the determination of contact $9$-dimensional filiform Lie algebras permits a quickly determination of the class of symplectic filiform $8$-dimensional Lie algebra. 
\medskip

From Proposition \ref{contact9}, we can highlight a model of $9$-dimensional filiform contact Lie algebra, that is a Lie algebra such as any $9$-dimensional filiform contact Lie algebra is isomorphic to a deformation of this model. We consider the Lie algebras $\g_{a_2,a_4,a_6}$ given by 
$$a_5=a_7=a_8=a_9=a_{10}=a_{11}=0$$
and
$$-3a_4^2+2a_6^2+2a_2a_6+a_4a_6=0, \ \ a_2a_4a_6 \neq 0.$$ Since $  a_6 \neq 0$, $a_2=\ds\frac{(a_4-a_6)(3a_4+2a_6)}{2a_6}$ and the condition $a_2a_4a_6 \neq 0$ is equivalent to $a_4a_6(a_4-a_6)(3a_4+2a_6)\neq 0.$
\begin{proposition}
Any $9$-dimensional filiform contact Lie algebra is isomorphic to a linear deformation of a Lie algebra $\g_{a_2,a_4,a_6}$.
\end{proposition}
\begin{corollary}
Any $8$-dimensional symplectic filiform Lie algebra is isomorphic to a linear symplectic deformation of a Lie algebra of the following family
$$
\left\{
\begin{array}{ l}
    \mu(X_0, X_i)=X_{i+1},  \ 1 \leq i \leq 6 ,\\
  \mu(X_1, X_5)= b_2 X_7\\
  \mu(X_2, X_4)=  b_4X_7,\\
\mu(X_1, X_4)= (b_2+b_4) X_6,\\
\mu(X_2, X_3)= b_4 X_6,\\
\mu(X_1, X_3)= (b_2+2b_4) X_5,\\
\mu(X_1, X_2)= (b_2+2b_4) X_4,\\
\end{array}
\right.
$$
with  $b_2=b_4=0$ or $b_4 ( b_2+b_4 )(2b_2-b_4)(b_2+2b_4)\neq 0.$
\end{corollary}

\subsection{The varieties $\mathcal{F}il_{10}$ and $\mathcal{F}il_{11}$}

\begin{proposition}Any $10$-dimensional filiform Lie algebra over $\K$ is given in a Vergne basis by 

\begin{equation}
\label{fili10}
\left\{
\begin{array}{ l}
    \mu(X_0, X_i)=X_{i+1},  \ 1 \leq i \leq 8 ,\\
    
     \mu(X_2, X_7)= a_1 X_9, \\
       \mu(X_1, X_7)= a_1 X_8+a_2X_9, \\
       
       \mu(X_3, X_6)= -a_1 X_9,\\
  \mu(X_2, X_6)= a_4X_9,\\
 \mu(X_1,X_6)=a_1X_7+(a_2+a_4)X_8+a_5X_9,\\
 
  \mu(X_4, X_5)= a_1 X_9,\\
  \mu(X_3, X_5)= a_7X_9,\\     
       \mu(X_2, X_5)= (a_4+a_7) X_8+a_8X_9,\\
  \mu(X_1, X_5)= a_1X_6+(a_2+2a_4+a_7) X_7+(a_5+a_8)X_8+a_9X_9,\\
  
  \mu(X_3, X_4)= a_7X_8+a_{10}X_9\\
  \mu(X_2, X_4)=  (a_4+2a_7) X_7+(a_8+a_{10})X_8+a_{11}X_9,\\
\mu(X_1, X_4)= a_1X_5+(a_2+3a_4+3a_7) X_6+(a_5+2a_8+a_{10})X_7+(a_9+a_{11})X_8+a_{12}X_9,\\

\mu(X_2, X_3)= (a_4+2a_7) X_6+(a_8+a_{10})X_7+a_{11}X_8+a_{13}X_9,\\
\mu(X_1, X_3)= a_1X_4+(a_2+4a_4+5a_7) X_5+(a_5+3a_8+2a_{10})X_6+(a_9+2a_{11})X_7+\\
\qquad \qquad \qquad(a_{12}+ a_{13})X_8+a_{14}X_9,\\
\mu(X_1, X_2)= a_1X_3+(a_2+4a_4+5a_7) X_4+(a_5+3a_8+2a_{10})X_5+(a_9+2a_{11})X_6+\\
\qquad \qquad \qquad(a_{12}+a_{13})X_7+ a_{14}X_8+a_{15}X_9,\\
\end{array}
\right.
\end{equation}
with the conditions
$$
\left\{\begin{array}{l}
     a_1(2a_2+7a_4+7a_7)=0,    \\
 3a_4^2 +3a_4a_7-2a_2a_7=0,\\
   a_1(2a_9+5a_{11})-2a_2a_{10}+a_4(7a_8-2a_{10})+a_7(-3a_5+2a_8-7a_{10})=0.  
\end{array}
\right.$$ 
\end{proposition} 
We denote the set of this multiplications by $Fil_{10}$.
If $2a_2+7a_4+7a_7=0$, then  
$$3a_4^2 +3a_4a_7-2a_2a_7=3a_4^2+10a_4a_7+7a_7^2=(a_4+a_7)(3a_4+7a_7).$$
We deduce
\begin{proposition}
The set $\mathcal{F}il_{10}$ of $10$-dimensional filiform Lie algebras is the union of the algebraic components
\begin{enumerate}
  \item $\mathcal{F}il_{10}(1)=\mathcal{O}(Fil_{10}(1))$ where $Fil_{10}(1)$ is the set of multiplication $\mu \in Fil_{10}$ satisfying
 $$ a_1=0, 3a_4^2 +3a_4a_7-2a_2a_7=0, -2a_2a_{10}+a_4(7a_8-2a_{10})+a_7(-3a_5+2a_8-7a_{10})=0.  $$
   \item $\mathcal{F}il_{10}(2)=\mathcal{O}(Fil_{10}(2))$ where $Fil_{10}(2)$ is the set of multiplication $\mu \in Fil_{10}$ satisfying
   $$ a_1, a_2=0, \ a_4=-a_7, \   a_1(2a_9+5a_{11})+a_4(3a_5+5a_8+5a_{10})=0,$$  
  \item $\mathcal{F}il_{10}(3)=\mathcal{O}(Fil_{10}(3))$ where $Fil_{10}(3)$ is the set of multiplication $\mu \in Fil_{10}$ satisfying
   $$ a_2=-2a_4, \ 3a_4=-7a_7, \   a_1(2a_9+5a_{11})+a_4\left( \frac{9}{7}a_5+\frac{43}{7}a_8+5a_{10}\right) =0.$$  
\end{enumerate}
\end{proposition}
\begin{proposition}Any $11$-dimensional filiform Lie algebra over $\K$ is given in a Vergne basis by 

\begin{equation}
\label{fili11}
\left\{
\begin{array}{ ll}
    \mu(X_0, X_i)=&X_{i+1},  \ 1 \leq i \leq 9 ,\\
      \mu(X_1, X_8)=& a_2X_{10}, \\

     \mu(X_2, X_7)= &a_4 X_{10}, \\
       \mu(X_1, X_7)= &(a_2+a_4) X_9+a_5X_{10}, \\
       
       \mu(X_3, X_6)= &a_7 X_{10},\\
  \mu(X_2, X_6)= &(a_4+a_7)X_9+a_8X_{10},\\
 \mu(X_1,X_6)=&(a_2+2a_4+a_7)X_8+(a_5+a_8)X_9+a_9X_{10},\\
 
  \mu(X_4, X_5)= &a_{10}X_{10},\\
  \mu(X_3, X_5)= &(a_7+a_{10})X_9+a_{11}X_{10},\\     
       \mu(X_2, X_5)=& (a_4+2a_7+a_{10}) X_8+(a_8+a_{11})X_9+a_{12}X_{10},\\
  \mu(X_1, X_5)=& (a_2+3a_4+3a_7+a_{10})X_7+(a_5+2a_8+a_{11}) X_8+(a_9+a_{12})X_9+a_{13}X_{10},\\
  
  \mu(X_3, X_4)=& (a_7+a_{10})X_8+a_{11}X_9+a_{14}X_{10},\\
  \mu(X_2, X_4)= & (a_4+3a_7+2a_{10}) X_7+(a_8+2a_{11})X_8+(a_{12}+a_{14})X_9+a_{15}X_{10},\\
\mu(X_1, X_4)=& (a_2+4a_4+6a_7+3a_{10}) X_6+(a_5+3a_8+3a_{11})X_7+(a_9+2a_{12}+a_{14})X_8\\
&+(a_{13}+a_{15})X_9+a_{16}X_{10},\\

\mu(X_2, X_3)= &(a_4+3a_7+2a_{10}) X_6+(a_8+2a_{11})X_7+(a_{12}+a_{14})X_8+a_{15}X_9+a_{17}X_{10},\\
\mu(X_1, X_3)=& (a_2+5a_4+9a_7+5a_{10}) X_5+(a_5+4a_8+5a_{11})X_6+(a_9+3a_{12}+2a_{14})X_7\\
&+(a_{13}+2a_{15})X_8+(a_{16}+a_{17})X_9+a_{18}X_{10},\\
\mu(X_1, X_2)=& (a_2+5a_4+9a_7+5a_{10}) X_4+(a_5+4a_8+5a_{11})X_5+(a_9+3a_{12}+2a_{14})X_6\\
&+(a_{13}+2a_{15})X_7+(a_{16}+a_{17})X_8+a_{18}X_{9}+a_{19}X_{10},\\
\end{array}
\right.
\end{equation}
with the conditions
$$
\left\{\begin{array}{l}
(J_1):  3z_2^2+3z_2z_3-2z_1z_3=0,    \\
  (J_2):  z_7(2z_1+2z_2+z_3)+z_3(3z_5+z_6) -7z_2z_6=0,\\
(J_3) : z_4(2z_1+7(z_2+z_3))-2z_3(2z_2+z_3)=0,\\
(J_4):  z_4(2z_8+5z_9)-z_{10}(2z_1+9z_2+12z_3 )-z_{7}(3z_5+7z_6-z_7)+4z_6^2  \\ \ \ \ \qquad \qquad \quad -2z_3(2z_8+7z_9)+8z_9(z_2+2z_3)=0. \\
\end{array}
\right.$$ 
where $z_1=a_2+a_4,z_2=a_4+a_7,z_3=a_7+a_{10},z_4=a_{10},z_5=a_5+a_{11},z_6=a_8+a_{11},z_7=a_{11},z_8=a_9+a_{12},z_9=a_{12}+a_{14}, z_{10}=a_{14}.$
\end{proposition} 
This system determines an irreducible component in $Fil_{11}$. This component is the Zariski closure of the open set whose elements are the filiform contact Lie algebra. 

\subsection{Contact and symplectic structures}
Let $\g$ be a $11$-dimensional filiform Lie algebra belonging to  $Fil_{11}$. Let $\{ \omega_1,\cdots,\omega_{10}\}$ the dual basis of the Vergne basis of $\g$. Assume that $\g$ is a contact Lie algebra. Then, from \cite{GKM}, the form $\omega_{10}$ is also a contact form. Then $\g$ is a contact algebra if and only if $\omega_{10}$ is a contact form in $\g_1$. We deduce

\begin{proposition}
A $11$-dimensional filiform Lie algebra. is a contact Lie algebra if and only if it is isomorphic to a Lie algebra  of $Fil_{11}$ with $a_2a_4a_7a_{10} \neq 0.$
\end{proposition}
We deduce that a model of contact $11$-dimensional filiform Lie algebra is given by the family $\g_{a_2,a_4,a_7,a_{10}}$ of Lie algebras of $Fil_{11}$ corresponding to
$$a_5=a_8=a_9=a_{11}=a_{12}=a_{13}=a_{14}=a_{15}=a_{16}=a_{17}=a_{18}=a_{19}=0$$
with the conditions,  
$$
\left\{\begin{array}{l}
 3z_2^2+3z_2z_3-2z_1z_3=0,   \\
z_4(2z_1+7z_2+7z_3)-2z_3(2z_2+z_3)=0,\\ 
  a_2a_4a_7a_{10}=z_4(z_3-z_4)(z_2-z_3+z_4)(z_1-z_2+z_3-z_4) \neq 0.
\end{array}
\right.$$ 
and $z_1=a_2+a_4,z_2=a_4+a_7,z_3=a_4+a_{10},z_4=a_{10}.$ Let us note that the algebraic variety determined by the two first equations is not reduced to $0$, for example the point $(a_2,a_4,a_7,a_{10})=(1,-1,1,-1)$ belongs to this algebraic set. Let us note also that the linear form $\omega_{10}$ where $\{\omega_i\}$ is the dual basis of $\{X_i\}$ satisfies
$$d\omega_{10}=-\omega _0 \wedge \omega_9 -a_2 \omega _1 \wedge \omega_8 -a_4 \omega _2 \wedge \omega_7-a_7 \omega _3 \wedge \omega_6-a_{10} \omega _4 \wedge \omega_5$$
and is a contact form.

\noindent Let $\mathcal{A}$ be the open set of $11$-dimensional filiform contact Lie algebra. Then this open set is the orbit of the family of Lie algebras $\g_{a_2,a_4,a_7,a_{10}}$  and satisfying $a_2a_4a_7a_{10} \neq 0$. Moreover 
$$\mathcal{F}il_{11}=\overline{\mathcal{A}}.$$

\begin{proposition}
None of the  $11$-dimensional filiform Lie algebras is rigid.
\end{proposition}
\pf The proof is similar to the $9$-dimensional case. The affine scheme being reduced, it is sufficient to prove that the dimension of the second space of cohomology $H^2_{CR}(\mu,\mu)$ of any Lie algebras of (\ref{fili11}) is not zero. It is clear that, if $\mu$ is a Lie algebra of (\ref{fili11}), the dimension of the $2$-cocycles $\psi$ such that $\mu +t \psi$ belongs to \ref{fili11} is of dimension $16$ parametrized by the $a_i.$ If $f$ is an endomorphism of $\K^{11}$, then putting $f(e_0)=\sum_{0}^{10} \alpha_iX_i$ and $f(e_1)=\sum_{1}^{10} \beta_iX_i$,  $\delta f (X_1,X_8)=v_2X_{10}$, $\delta f (X_2,X_7)=v_4X_{10}$, $\delta f (X_3,X_6)=v_7X_{10}$, $\delta f (X_4,X_5)=v_{10}X_{10}$, we have
$$
\left\{
\begin{array}{l}
  v_2=  a_2(\beta_1-2\alpha_0)  \\
   v_4=   a_4(\beta_1-2\alpha_0) \\
     v_7=  a_7(\beta_1-2\alpha_0)  \\
       v_{10}=  a_{10}(\beta_1-2\alpha_0)  \\  
\end{array}
\right.
$$
We deduce that the dimension of the space of deformations is greater or equal to $3$.

\medskip

\noindent{\bf Consequence: Determination of the symplectic $10$-dimensional filiform Lie algebras.}  From the Proposition \ref{consym} we deduce:

\begin{proposition}
Any $10$-dimensional symplectic filiform Lie algebra is isomorphic to
$$
\left\{
\begin{array}{ ll}
    \mu(X_0, X_i)=&X_{i+1},  \ 1 \leq i \leq 8 ,\\

       \mu(X_1, X_7)= &(a_2+a_4) X_9, \\
       
  \mu(X_2, X_6)= &(a_4+a_7)X_9,\\
 \mu(X_1,X_6)=&(a_2+2a_4+a_7)X_8+(a_5+a_8)X_9,\\

  \mu(X_3, X_5)= &(a_7+a_{10})X_9,\\     
       \mu(X_2, X_5)=& (a_4+2a_7+a_{10}) X_8+(a_8+a_{11})X_9,\\
  \mu(X_1, X_5)=& (a_2+3a_4+3a_7+a_{10})X_7+(a_5+2a_8+a_{11}) X_8+(a_9+a_{12})X_9,\\
  
  \mu(X_3, X_4)=& (a_7+a_{10})X_8+a_{11}X_9,\\
  \mu(X_2, X_4)= & (a_4+3a_7+2a_{10}) X_7+(a_8+a_{11})X_8+(a_{12}+a_{14})X_9,\\
\mu(X_1, X_4)=& (a_2+4a_4+6a_7+3a_{10}) X_6+(a_5+3a_8+a_{11})X_7+(a_9+a_{12}+a_{14})X_8\\
&+(a_{13}+a_{15})X_9,\\

\mu(X_2, X_3)= &(a_4+3a_7+2a_{10}) X_6+(a_8+2a_{11})X_7+(a_{12}+a_{14})X_8+a_{15}X_9,\\
\mu(X_1, X_3)=& (a_2+5a_4+9a_7+5a_{10}) X_5+(a_5+4a_8+5a_{11})X_6+(a_9+3a_{12}+2a_{14})X_7\\
&+(a_{13}+2a_{15})X_8+(a_{16}+a_{17})X_9,\\
\mu(X_1, X_2)=& (a_2+5a_4+9a_7+5a_{10}) X_4+(a_5+4a_8+5a_{11})X_5+(a_9+3a_{12}+2a_{14})X_6\\
&+(a_{13}+2a_{15})X_7+(a_{16}+a_{17})X_8+a_{18}X_{9}\\
\end{array}
\right.
$$
with $a_2a_4a_7a_{10}\neq 0$ and if $(z_1=a_2+a_4,z_2=a_4+a_7,z_3=a_7+a_{10},z_4=a_{10},z_5=a_5+a_{11},z_6=a_8+a_{11},z_7=a_{11},$
$$
\left\{\begin{array}{l}
 3z_2^2+3z_2z_3-2z_1z_3=0,    \\
 z_7(2z_1+2z_2+z_3)+z_3(3z_5+z_6) -7z_2z_6=0.\\
\end{array}
\right.$$ 
\end{proposition}
We deduce, from the definition of a symplectic model:
\begin{corollary}
The symplectic models of  $10$-dimensional filiform symplectic Lie algebras are the Lie algebra  corresponding to
$$a_5=a_8=a_{11}=a_9=a_{12}=a_{14}=a_{12}=a_{14}=a_{15}=a_{16}=a_{17}=a_{18}=0.$$
\end{corollary}

\section{Contact and symplectic filiform Lie algebras}

\subsection{$(2p+1)$-dimensional contact filiform Lie algebras}

Let $\{X_0,\cdots,X_{2n}\}$ be a Vergne basis of a $(2p+1)$-dimensional filiform Lie algebra $\g$ and let us denote by $C_{i,j}^k$ the constant structures related to this basis. We have seen in dimension $11$ or smaller, but this remains trivially true for greater dimension, that the structure constants of $\g$ related to this basis are linear combinations of the $(p-1)^2$ structure constants $a_{i,j}=C_{i,j}^{2p}$ for   $1 \leq i < j \leq 2p-1-i.$ Using the same notation as the previous section, we deduce that $Fil_{2p+1}$ is an algebrac variety embedded in $\K^{(p-1)^2}$.
In fact, all the other structure constants are defined by the linear equation
$$[X_0,[X_i,X_j]]=[X_{i+1},X_j]+[X_i,X_{j+1}]$$
as soon as $j+1 \leq 2p$. More precisely, we have
$$
\left\{
\begin{array}{l}
\medskip
     C_{2p-1-j-k,j}^{2p-k}=\sum_j \alpha_{j}a_{2p-1-j,j}, \ \ 2j > 2p-1-k, \ \ k=1,\cdots,2p-4    \\
     \medskip
     C_{2p-2-j-k,j}^{2p-k}=\sum_j \beta_{j}a_{2p-2-j,j}, \ \ 2j > 2p-2-k, \ \ k=1,\cdots,2p-5     \\
     \cdots \\
     C_{2p-(2p-3)-j-k,j}^{2p-k}=C^{2p-k}_{3-j-k,j}=a_{1,2},     \\
     \end{array}
\right.
$$
The coefficients $\alpha_j,\beta_j$ are described in \cite{Goze-Khakim}. Let $\{\omega_0,\omega_1,\cdots,\omega_{2p}\}$ be the dual basis. All the Jacobi conditions are given writing that $d(d\omega_{2p})=0.$ This gives 
\begin{equation}
\label{Jac}
(p-3)^2+(p-4)(p-5)+(p-6)^2+ \cdots +\varepsilon
\end{equation}
equations where $\varepsilon=2$ if $p \equiv 0 \pmod 3$, $\varepsilon=1$ if $p \equiv 1 \pmod 3$ and $\varepsilon=2^2$ if $p\equiv 2 \pmod 3$.  This shows that, as soon as the dimension exceeds $19$ the number of  polynomial equations is greater than the number of parameters $a_{i,j}$.   

If $\g$ admits a contact form, then from \cite{GKM} the form $\omega_{2p}$ is also a contact form. We deduce

\begin{proposition}
 A $(2p+1)$-dimensional filiform Lie algebra admits a contact form if and only if the structure constants related to the Vergne basis satisfyy:
 $$a_{1,2p-2}\cdot a_{2,2p-3}\cdots a_{i,2p-1-i}\cdots a_{p-1,p} \neq 0.$$
 \end{proposition}
Since any deformation of a contact Lie algebra is also a contact Lie algebra, we deduce that the set of  $(2p+1)$-dimensional filiform contact Lie algebra is a Zariski open set in $\mathcal{F}il_{2p+1}$.  Let us consider the family $\mathcal{A}$ contained in this open set and corresponding to Lie algebras whose structure constants satisfy $a_{i,j}=0$ except $a_{1,2p-2}, a_{2,2p-3}, \cdots,  a_{i,2p-1-i},\cdots,  a_{p-1,p}$ which are supposed different to $ 0.$ It is clear that any contact filiform $(2p+1)$-dimensional Lie algebra is a deformation of a Lie algebra of this family. We remark that this family is parametrized by the $(p-1)$ constant structures $a_{1,2p-2}, a_{2,2p-3}, a_{i,2p-1-i},  a_{p-1,p}$ but the system of polynomial equations deduced from the Jacobi conditions, which is a consequence of $d(d\omega_{2p})=0$ is composed, when $p$ is greater than $7$, of a number of equations greater than $p-1$. This number depend to $p \ mod(3)$. For example, if $p=3k+1$, we have $3k$ parameters and $3k^2-3k+1$ polynomial equations and $3k^2-3k+1 >3k$ as soon as $k \geq 2$.

Let us note also that $\mathcal{A}$ is not reduced to $0$. In fact the Lie algebra corresponding to
\begin{equation}
\label{g0}
a_{1,2p-2}=1,a_{2,2p-3}=-1,\cdots, a_{i,2p-1-i}=(-1)^{i+1},\cdots, a_{p-1,p} =(-1)^p
\end{equation}
belongs to this family and more generally, if
$$a_{1,2p-2}=\lambda,a_{2,2p-3}=-\lambda,\cdots, a_{i,2p-1-i}=(-1)^{i+1}\lambda,\cdots, a_{p-1,p} =(-1)^p\lambda$$
then the corresponding Lie algebras are in $\mathcal{A}$. This implies that $\dim \mathcal{A} \geq 1.$ Let us denote by $\g_0$ the Lie algebra of $\mathcal{A}$ corresponding to (\ref{g0}) and let us compute the space of deformations. As we have seen in previous work (\cite{GRKegel}), it is sufficient to compute the cocycles of $\g_0$ which preserves the Vergne's basis. For example, in case of dimension $9$, from (\ref{fili9}), the space of such cocycles is of dimension $8$. Let us compute the space of coboundaries. It is generated by the $\delta f$ satisfying $\delta f(X_0,Y)=0$ and $f \in gl(9,\K)$. This last condition implies that $f$ is determined when we know $f(X_0)=\alpha_0 X_0 +\cdots+\alpha_8 X_8$ and $f(X_1)=\beta_1 X_1 +\cdots + \beta_8 X_8.$ We obtain
$$\delta f(X_1,X_6)=(-2\alpha_0+ \beta_1)X_8= -\delta f(X_2,X_5)=\delta f(X_3X_4)$$
and
$$\delta f(X_1,X_4)=-\delta f(X_2,X_3)=-2\beta_3X_8, \ \delta f(X_1,X_2)=-2\beta_5 X_8,$$
if not $\delta f(X_i,X_j)=0.$ We deduce that the space of deformations is of dimension $4$. In dimension $2p$, the space of cocycles which preserves the Vergne's basis is embedded in a vector space of dimension $(p-1)^2$ parametrized by the structure constants 
$$\{a_{1,2p-2},\cdots,a_{1,2},a_{2,p-3},\cdots,a_{2,3},\cdots,a_{p-2,p},a_{p-2,p-1},a_{p-1,p}\}$$
that is $a_{i,j}$ with $1 \leq i < j \leq 2p-2$, $3 \leq i+j \leq 2p-1$ and with $a_{i,j}=C_{i,j}^{2p}.$
  If $f$ is a linear endomorphism of $gl(2p,\K)$, then $\delta f(X_0,X_i)=0$ implies that $f(X_i)$ is determined for $2 \leq i \leq 2p$ by $f(X_0)=\alpha_0 X_0 +\cdots+\alpha_{2p} X_{2p}$ and $f(X_1)=\beta_1 X_1 +\cdots + \beta_{2p} X_{2p}.$ This implies
$$\delta f(X_1,X_{2p-2})=(-2\alpha_0+ \beta_1)X_{2p}= -\delta f(X_2,X_{2p-3})=\cdots =(-1)^p\delta f(X_{p-1}X_{p}),$$
and the  other non zero $\delta f(X_i,X_j)$ are
$$
\left\{
\begin{array}{l }
   \delta f(X_1,X_{2i}) =2\beta_{2p-2i-1}X_{2p}, \ 1 \leq i \leq p-2    \\
    \delta f(X_2,X_{2i+1}) =-2\beta_{2p-2i-3}X_{2p}, \ 1 \leq i \leq p-3  \\
     \delta f(X_3,X_{2i}) =2\beta_{2p-2i-5}X_{2p} , \ 2 \leq i \leq p-2   \\
     \cdots \\
     \delta f(X_{p-2},X_{p-1})=\beta_3 X_{2p}.
\end{array}
\right.
$$

We remark also than the parameters $a_{1,2},a_{1,3},a_{2,3},a_{1,4},a_{2,4},a_{1,5}$ do not appear in the polynomial Jacobi equations because the forms $\omega_i \wedge \omega_j$ which appear in the expression of $d\omega_{2^p}$ are closed for $(i,j) \in \{(1,2),(1,3),(2,3),(1,4),(2,4),(1,5)\}.$  From the previous computations of the coboundaries, we have in particular
$$
\left\{
\begin{array}{l } 
\delta f(X_1,X_{2}) =2\beta_{2p-3}X_{2p}, \    \delta f(X_1,X_{4}) =2\beta_{2p-5}X_{2p} , \ \delta f(X_1,X_{3}) =\delta f(X_1,X_{5}) =0,\\
\delta f(X_2,X_{3}) =2\beta_{2p-5}X_{2p}, \    \delta f(X_2,X_{4}) =0.
\end{array}
\right.
$$
Then we can consider that the parameters $a_{1,2},a_{1,4}$ are orbital parameters and $a_{1,3},a_{2,3},a_{2,4},$ $a_{1,5}$ are parameters of non trivial deformations. To end this work, we compute the space of deformation of $\g_0$.
It remains to compute the dimension of the space of cocycles.  We have seen that it is embedded in a vector space of dimension $(p-1)^2$ and the number of polynomial equations given by the Jacobi relations was (\ref{Jac}). But the affine scheme associated to this polynomial equations is not reduced. We can find relations between these equations from the following remark which we illustrate in dimension $11$. In this dimension the Jacobi polynomial equations is constituted of $4$ equations corresponding to $J(X_i,X_j,X_k)=0$ ($J$ for the Jacobi condition related to triple $(X_i,X_j,X_k)$. To simplify we denote  by $(i,j,k)$ the polynomial $J(X_i,X_j,X_k)$ and let $ p=i+j+k$ be its weight. In this case we have $4$ parameters $(a_{1,8}=a_2,a_{2,7}=a_4,a_{3,6}=a_7,a_{4,5}=a_{10})$ and $4$ equation corresponding to $(p=6, (i,j,k)=(1,2,3))$, $(p=7, (i,j,k)=(1,2,4))$, $(p=8, (i,j,k)=(1,2,5), (1,3,4)).$ But we have
$$\begin{array}{l}
     [X_0,(1,2,3)]=(1,2,4)    \\
      \lbrack X_0,(1,2,4)]=(1,3,4) +(1,2,5).  
\end{array}
$$
Thus the system of Jacobi equations can be reduced to the system
$$(1,2,3)=0, (1,3,4)=0$$
and the corresponding affine scheme is reduced. Then we have $4$ parameters which have to satisfy $2$ independent relations. The space of parameters is then of dimension $2$. Let us come back to the general model $\g_0$. The linear space of parameters is of dimension $(p-1)$ and it is generated by the structure constants
$$a_{1,2p-2},a_{2,2p-3},\cdots,a_{p-2,p+1},a_{p-1,p}. $$
The  weights take its values in $(6,7,\cdots,2p-2)$ and concern the Jacobi equation:
$$(1,2,3); (1,2,4); (1,2,5),(1,3,4);\cdots;(1,p-2,p-1),(1,p-3,p),\cdots,(1,2,2p-5),$$
$$(2,p-1,p+1),
\cdots,(2,3,2p-7),\cdots,\}$$
the last term in this ordered sequence depends to $(p \  \mod 3)$, more precisely, if $2p-2  \equiv 2 \pmod 3$, then the last term is $(k-1,k+1,k+2)$ with $2p-2=3k+2,$ if $2p-2  \equiv 1 \pmod 3$, then the last term is $(k-1,k,k+2)$, and if $2p-2  \equiv 0 \pmod 3$, then the last term is $(k-1,k,k+2)$. We have seen (\ref{Jac}) that the number of Jacobi polynomial equations is 
$$(p-3)^2+(p-4)(p-5)+(p-6)^2+ \cdots +\varepsilon$$
where $\varepsilon=2$ if $p \equiv 0 \pmod 3$, $\varepsilon=1$ if $p \equiv 1 \pmod 3$ and $\varepsilon=2^2$ if $p\equiv 2 \pmod 3$.
This scheme is not reduced. To reduce it we consider the relations
$$[X_0,(i,j,k)]=(i+1,j,k)+(i,j+1,k)+(i,j,k+1).$$
Putting $2p-2=3m+r$ with $0 \leq r \leq 2$, we can write the reduced number $N_r$ of equations:

\begin{itemize}
  \item If $m=2h$  and $r=0$, then $N_r=3h^2-3h+1$,
  \item If $m=2h+1$ and $r=1$, then $N_r=3h^2+3h$,
  \item If $m=2h$ and $r=2$, then $N_r=3h^2-h$.
\end{itemize}
We can see than, as soon as $n \geq 14$, that is $p=7,m=4,r=2$ then the number of parameters is $6$ and $N_r=10.$  In the same way, we can reduce this new polynomial system using the identity
$$[X_1,(i,j,k)]=((i+2,j,k)C_{1,i}^{i+2}+(i,j+2,k)C_{1,j}^{j+2}+(i,j,k+2)C_{1,k}^{k+2})X_{i+j+k+2}.$$
which is a direct consequence of the natural grading of $\g_0$. We deduce in particular
$$[X_1,(1,2,3)]=(-(1,3,4)C_{1,2}^{4}+(1,2,5)C_{1,3}^{5})X_{10}.$$
  To end this section, we can look the case $n=14$, this case corresponding to $N_r > p-1$. We have $6$ coefficients and $7$ relations after the reduction of the first type. We can choose as generating relations, the relation $(1,2,i)$ for $i=3,5,6,7,8,9$ and $(3,4,5)$. The relation of second type concerning these equations are, where $(1,(i,j,k)$ is the coefficient of $[X_1,(i,j,k)]$, 
$$
\left\{
\begin{array}{l}
  1(1,2,3)=(1,2,5)C_{1,3}^5 ,    \\
  1(1,2,4)=(1,2,6)C_{1,4}^6,     \\
  1(1,2,5)=(1,4,5)C_{1,2}^4+(1,2,7)C_{1,5}^7,     \\
    1(1,2,6)=(1,4,6)C_{1,2}^4+(1,2,8)C_{1,6}^8,     \\
  1(1,2,7)=(1,4,7)C_{1,2}^4+(1,2,9)C_{1,7}^9,     \\
    1(2,3,5)=-(3,4,5)C_{1,2}^4+(2,3,7)C_{1,5}^7,     \\
  \end{array}
  \right.
  $$
  If $C_{1,2}^4C_{1,3}^5 C_{1,4}^6C_{1,5}^7C_{1,6}^8C_{1,7}^9 \neq 0$, then the Jacobi polynomial system is reduced only to one equation $(1,2,3)$. In this open set, the space of parameters of deformations of the models is of dimension greater or equal to $5$. 
  
  \subsection{Filiform symplectic algebras}
  
  From the previous study, we have that the $(2p)$-dimensional symplectic filiform Lie algebras are isomorphic to a quotient of a contact $(2p+1)$-dimensional filiform Lie algebra $\g_{2p+1}$  by its center $\mathbb{K}\{X_{2p}\}$. Then its can be written with the structure constants of $\g_{2p+1}$ with the condition $a_{1,2p-2}a_{2,2p-4}\cdots a_{p-1,p} \neq 0.$


\end{document}